\documentclass[12pt,a4paper,english]{article}
\usepackage[T1]{fontenc}
\usepackage[utf8]{inputenc}
\usepackage{babel}
\usepackage{amsthm}
\usepackage{amsmath}
\usepackage{amssymb}
\usepackage{amsfonts}
\usepackage{mathpazo}
\usepackage{newtxtext}
\usepackage[normalem]{ulem}
\usepackage[all]{xy}
\usepackage{tikz}
\usetikzlibrary{matrix,arrows,decorations.pathmorphing,positioning}
\usepackage{hyperref}
\usepackage{fundamentalgroupaffineline}


\begin{document}

\centerline{ } 
\centerline{\LARGE{The Geometric Fundamental Group}} 
\vspace{4pt} 
\centerline{\LARGE{of the Affine Line over a Finite Field}} 
\vspace{25pt}
\centerline{\large{Henrik Russell}} 
\vspace{10pt} 
\centerline{\large{\today}} 
\vspace{25pt}

\begin{abstract} 
The affine line $\Afn^1_{\finfld}$ 
and the punctured line $\Afn^1_{\finfld} \setminus \st{0}$ 
over a finite field $\finfld$ 
are taken as benchmarks for the problem of describing 
(non-commutative) geometric \'etale fundamental groups. 
To this end, using a reformulation of Tannaka duality 
we construct for a projective variety $X$ 
a (non-commutative) universal affine pro-algebraic group $\LalbG_X$, 
such that for any given affine subvariety $U$ of $X$ 
any finite and \'etale Galois covering of $U$ over $\finfld$ 
is a pull-back of a Galois covering of a quotient $\LalbGs{X}{U}$ of $\LalbG_X$. 
Then the geometric fundamental group $\fundGgeoF{U}{\ol{\finfld}}$ 
is a completion of $\LalbGs{X}{U}\pn{\ol{\finfld}}$. 
We obtain explicit descriptions of the universal affine groups $\LalbGs{X}{U}$ 
for $U =$ $\Afn^1_{\finfld}$ and $\Afn^1_{\finfld} \setminus \st{0}$. 
\end{abstract}

\setcounter{tocdepth}{1}
\tableofcontents{} 

\setcounter{section}{-1}

\newpage

\section{Introduction} 
The \'etale fundamental group is the algebraic analogue of the topological fundamental group. 
For schemes corresponding to analytic manifolds 
the \'etale fundamental group is the completion of the topological fundamental group 
of the corresponding analytic manifold. 
In most other cases the shape of this group is an open question, 
e.g.\ for (geometrically non-trivial) schemes over a finite field. 
In other words: 
we have no idea what this group might look like in these cases 
it was initially designed for, 
not even for the simplest example one can think of, the affine line over a finite field. 
The aim of this paper is to give an idea of how this problem might be attacked. 

\bThm[Serre, Raynaud, Harbater] 
\label{quasi-p-group}
Let $\clfld$ be an algebraically closed field of characteristic $p > 0$. 
A finite group $\Gam$ appears as the Galois group 
of an unramified covering of the affine line \,$\Afn^1_{\clfld}$\, 
if and only if \,$\Gam$\, is a \emph{quasi-$p$-group} \linebreak 
(i.e., a group that is generated by its $p$-subgroups). 
\eThm 

This is a consequence of Abhyankar's conjecture for affine curves, 
which has been proven by Serre, Raynaud and Harbater 
(see \cite[Sec.\ 3]{HOPS}). 
This insight 
is so impressive that one might shy away from asking 
what the \'etale fundamental group of the affine line might look like as a whole. 
But looking at a large group in its entirety 
and examining the (possibly enormous) crowd of quotient groups 
are different tasks (cf.\ \cite[Introduction]{Ku}). 
When describing the motion of a body, it is easier to look at the barycentre 
than at each individual atom. 

The \emph{abelianized} fundamental group (classifying abelian coverings) 
can be characterized 
for curves by means of the generalized Jacobian of Rosenlicht-Serre, 
and for higher dimensional varieties by means of generalized Albanese varieties 
from \cite{Ru13}. 
The basic idea is to present Galois coverings (in the sense of rational maps) 
of a projective variety $X$ 
as pull-backs of isogenies to a universal pro-algebraic group $\Guniv_X$. 
The group of automorphisms of a Galois covering is then given by the kernel of an isogeny. 
Moreover, when working over a finite field 
there is an isogeny that is \emph{universal} in the sense that 
every isogeny to the universal group is a quotient of this universal isogeny. 
This universal isogeny is given by the Artin-Schreier map $\asm$ 
(= Frobenius minus identity), 
and the abelianized fundamental group is given by the kernel of $\asm$, 
which is the group of $k$-valued points of $\Guniv_X$ 
(cf.\ \cite[VI, No.\ 6--12, in part.\ No.\ 11, Prop.\ 9]{S_GrpAlg} and 
\cite[Thm.\ 8.8ff]{Ru_CFT}).

Is it possible to use a similar approach in the \emph{non-commutative} case? 
\marginpar{Sec.~\ref{sec:GaloisCov-via-Isogenies}}
For the purpose of Galois coverings, it is sufficient to work with 
a universal \emph{affine} group $\Lalb_X$ (see \cite[VI, No.~8]{S_GrpAlg}). 
Then the central problem is the construction \marginpar{Sec.~\ref{sec:Universal-Group}} 
of a \emph{non-commutative} universal affine pro-algebraic group $\LalbG_X$, 
i.e., a group satisfying the universal factorization property 
even for non-commutative groups. 
(This is accomplished in Theorems \ref{Exist_univAffObj} and \ref{univ_affObject}.) 

A main ingredient in the approach of \cite{Ru13} for a universal group 
is Cartier duality between affine groups and formal groups. 
This translates the group law of an affine group 
to the multiplication of the algebra of the corresponding formal group. 
The duality relates the affinization $\Lalb_X$ of the universal group $\Guniv_X$
to a formal subgroup $\HDivf_X$ 
(= principal divisors) 
of the sheaf of relative Cartier divisors $\Divf_X$ 
(cf.\ \cite{Ru_abFund}). 
Now if the universal affine group $\LalbG_X$ is supposed to be non-commutative, 
this means that we have to extend $\HDivf_X$ 
to a larger object with a non-commutative ``algebra multiplication''. 

The solution to this problem is straightforward, 
\marginpar{Sec.~\ref{sec:TannakaDuality}} 
we need to replace Cartier duality by a suitable reformulation of Tannaka duality: 
this is an anti-equivalence between affine groups and so called 
\emph{representation functors} (Theorem \ref{Tannaka duality}). 
The functor $\HDivf_X$ 
gives rise to a functor $\Autf_{X,\bpt}$ 
\marginpar{Pbl.~\ref{pbl: group fctr of rel Cartier div}} 
of based rational maps from $X$ to linear groups, 
and the subfunctor $\HDivsf{X}{S}$ of those divisors 
with support in a closed subset $S \subset X$ 
is replaced by the subfunctor $\Autfs{X,\bpt}{\,S}$ 
\marginpar{Sec.~\ref{sec:CategoryAut}} 
of those rational maps that are regular outside $S$ 
(see Definition \ref{group of based maps}). 
The dual $\LalbGs{X}{S}$ of $\Autfs{X,\bpt}{\,S}$ is characterized by 
a universal mapping property 
for based morphisms from $X \setminus S$ to affine algebraic groups. 
\marginpar{Sec.~\ref{sec:Recover from Repr}} 
We may recover $\LalbGs{X}{S}$ from $\Autfs{X,\bpt}{\,S}$ 
by means of Tannaka's Theorem (derived from \cite{Tannaka}), 
see \Point \ref{Repr from CartierDual}. 
The descriptions that we found are explicit enough 
to allow for the computation of concrete examples: 
this is done for $X \setminus S = \Afn^1$ and $\Afn^1 \setminus \st{0}$, 
see Examples \ref{HDivf(gl) for A^1}, \ref{HDivf(gl) for Gm}, 
\ref{L from HDiv for A^1} and \ref{L from HDiv for Gm}. 

\bigskip
Let $X$ be a projective curve and $U = X \setminus S$ an affine subvariety 
over a finite field $\fld$.
The group classifying geometric finite \'etale abelian coverings 
\marginpar{Sec.~\ref{sec:FundGroup}} 
of $U$ that are defined over $\fld$ 
is given by the $\fld$-valued points of the abelian universal affine group for $U$: 
\marginpar{Thm.~\ref{Abelian Coverings}} 
\[ \fundGabO{U} = \Lalbs{X}{S}\pn{\fld} 
\laurin 
\]
This is possible, since the Artin-Schreier map 
(the commutative version of the Lang map $\asm: x \lmt x^{-1} \Frob(x)$, 
where $\Frob$ is the Frobenius) is a homomorphism of groups, 
providing a universal isogeny 
(cf.\ \cite[VI, No.\ 6, Prop.\ 6]{S_GrpAlg}). 
In the non-commutative setting the Lang map is not necessarily a homomorphism. 
The group that classifies geometric finite and \'etale Galois coverings over $\fld$ 
is then 
\marginpar{Thm.~\ref{FundGroup-of-Affine}} 
\[ \fundGO{U} = \ProSub\bigpn{\LalbGs{X}{S}\pn{\fld}}
\] 
(proven here only for $U = \Afn^1$ and $\Afn^1 \setminus \st{0}$), 
where $\ProSub$ is a completion assigning to a pro-finite group $\Gam$ 
\marginpar{No.~\ref{pro-sub-completion}} 
the limit of the inverse system of subsets of quotients of $\Gam$. 
In other words: the set of finite quotients of \,$\fundGO{U}$\, 
is given by the set of subgroups of finite quotients of \,$\LalbGs{X}{S}\pn{\fld}$. 

In the following, let $\fld$ be a finite field 
and  $\clfld$ an algebraic closure. 
Puzzling together the pieces we have received so far, 
we can prove at the end of Section \ref{sec:FundGroup}: 

\bThm
\label{FundGroup-of-A1}
The geometric \'etale fundamental group of the affine line $\Afn^1$ 
over \,$\fld = \bF_p$ 
is generated by matrix-valued points of the big Witt vectors: 
\marginpar{Sec.~\ref{sec:FundGroup}} 
\[ \fundGgeo{\Afn^1} 
= \ProSub\biggpn{\varprojlim_n \Bigpn{1 + s\,\gl_n\pn{\clfld}[[s]]}^{\tms} } 
\]
where 
$s$ is a local parameter at \,$\infty \in \Prj^1 \setminus \Afn^1$ and 
\,$\gl_n = \ul{\Lie}\pn{\Gl_n} = \Mat_{n \tms n}$ is the algebra of $(n \tms n)$-matrices. 
The structure of an inverse system is induced by 
restriction to the upper left submatrix \,$\gl_m \la \gl_n$\, for $m < n$. 
\eThm 

A quick plausibility check is as follows. 
First, the abelian quotient $\fundGabgeo{\Afn^1}$ of $\fundGgeo{\Afn^1}$ 
is the pro-sub-completion of the $\clfld$-points of 
the group of (commutative) big Witt vectors: 
over $\fld$, the group $\fundGab\bigpn{\Afn^1_{\fld}}^0$\, is given by 
the $\fld$-valued points of the big Witt vectors 
\[ \fundGab\bigpn{\Afn^1_{\fld}}^0 = \Bigpn{1 + s\,\fld[[s]]}^{\tms} 
\]
($\fld$-valued points of generalized Jacobians 
as described in \cite[V, \S3, No.\ 14]{S_GrpAlg}, 
cf.\ \cite[Thm.\ 2.3]{Ru_abFund}), 
and $\fundGabgeo{\Afn^1}$ is the limit of an inverse system that contains the 
$\fundGab\bigpn{\Afn^1_{\erwfld}}^0$ as quotients 
for all finite extensions $\erwfld | \,\fld$. 
Second, 
every finite group is a subgroup of a matrix group, 
therefore the non-commutativity can be expressed in terms of matrix groups. 

Likewise, we can show at the end of Section \ref{sec:FundGroup}: 
\bThm
\label{FundGroup-of-Gm}
The geometric \'etale fundamental group of the punctured affine line 
$\Afn^1 \setminus \st{0}$ over a finite field \,$\fld = \bF_p$ is given by 
\marginpar{Sec.~\ref{sec:FundGroup}} 
\[ \fundG\bigpn{\Afn^1_{\clfld} \setminus \st{0}} 
= \ProSub\biggpn{\Gm\pn{\clfld} \tms 
	\varprojlim_n \Bigpn{1 + s\,\gl_n\pn{\clfld}[[s]]}^{\tms} } 
\laurin 
\]
\eThm 

Because of the intriguing fact that the candidate of Theorem \ref{FundGroup-of-A1} 
admits a quotient map to every finite \emph{quasi-$p$-group} 
due to Theorem \ref{quasi-p-group}, 
\marginpar{Sec.~\ref{sec:Quotient Maps}} 
we will illustrate in \Point \ref{quotient map} how to construct these. 
We will explain in \Point \ref{quotient map}, Step 2, 
why $\fundGabO{U}$, in contrast to $\fundGO{U}$, 
does \emph{not} need a pro-sub-completion: 
for an \emph{abelian} algebraic group $G$ over a finite field $\fld$, 
every subgroup of $G(\fld)$ is isomorphic to a quotient group of $G(\fld)$. 

\bigskip 
The relation between \'etale fundamental group and Tannakian fundamental group 
has already been a topic in other situations (see e.g.\ \cite{EsnaultPhung_fundGroupoid} or \cite{Deninger_fundGroup-TopolSpaces}), 
our approach however gives a different viewpoint to this 
(Problem \ref{pbl: group fctr of rel Cartier div}), 
\marginpar{Sec.~\ref{sec:Open Problems}} 
as well as provides a link to non-commutative class field theory 
(Problem  \ref{genl FundGroup of Affine}) 
and anabelian geometry (Problem \ref{Recover CoordRing from FundGroup}) 
for varieties over a finite field.

\vspace{5mm}
\noindent 
\textbf{Acknowledgements.} 

\noindent 
I thank David Harbater for 
helpful discussions.

\subsection{Notations and Conventions} 
\label{Terminology}

\subsubsection{Schemes} 
If $X$ is a scheme, then $\sK_X$ denotes the sheaf of total quotient rings 
of the structure sheaf $\sO_{X}$. 
If $X$ is irreducible, then $K_{X}$ denotes its function field. 
The coordinate ring of $X$ is denoted by $\sO(X) := \Gam\pn{X, \sO_X}$. 
A variety is an irreducible reduced separated scheme of finite type over a field. 
For any $k$-vector space $V$ we let 
\,$\ul{V} := \Spec\bigpn{\Sym\pn{V^{\vee}}}$ 
be the associated vectorial group \,$R \lmt V \tens R = \Homk\pn{V^{\vee}, R}$. 
The category of $k$-vector spaces is denoted by \,$\Veck$.

\subsubsection{Pro-Sub-Completion} 
\label{pro-sub-completion} 
Let \,$\Gam = \varprojlim \Gam_i$\, be the inverse limit of an inverse system of groups 
\,$\st{\Gam_i}_{i \in I}$\, over a partially ordered index set $I$, 
with transition maps \,$\gam_{ij}: \Gam_j \ra \Gam_i$\, for \,$i < j$. 
Let \,$\ProSub\pn{\st{\Gam_i}_{i \in I}}$\, be the inverse system given by 
the family \,$\st{\Sig \;|\; \Sig \subset \Gam_i, i \in I}$\, 
of subgroups of \,$\st{\Gam_i}_{i \in I}$, 
with induced maps \,$\sig_{ij} := \gam_{ij}|_{\Sig}: \Sig \lra \gam_{ij}\pn{\Sig}$.  
Then \,$\ProSub(\Gam) := \varprojlim \ProSub\pn{\st{\Gam_i}_{i \in I}}$\, 
is the inverse limit of the inverse system \,$\ProSub\pn{\st{\Gam_i}_{i \in I}}$. 
Note that this limit is not necessarily directed.


\subsubsection{Galois Coverings and Classifying Groups}
A rational map $\morphism: Y \dra X$ of varieties over $\fld$ 
is called a \emph{Galois covering} 
if it induces a Galois extension of function fields $K_Y|K_X$. 
It is called \emph{geometric}, if $\fld$ is algebraically closed in $K_Y$, 
or equivalently if \,$\Gal\pn{K_Y|K_X} \cong \Gal\pn{K_Y \clfld \big| K_X \clfld}$, 
i.e.\ the covering does not arise from an extension of the base field.  
The based geometric Galois coverings of $(X,\bpt)$ 
that are defined over $\fld$ 
are classified by the following group: 
\begin{align*}
	\bGalO{X,\bpt} 
	:=& \varprojlim_{\substack{(Y,\bptt) \dra (X,\bpt) \\ \textrm{ geometric Galois}/\fld}} \Gal\pn{K_Y|K_X} 
	\laurin 
\end{align*}
We define the group \,$\bGalgeo{X,\bpt}$\, as an inverse limit as above, 
but we drop the restriction that the Galois coverings are defined over $\fld$, 
i.e., the field of definition ranges over finite field extensions \,$\erwfld | \,\fld$\, 
in an algebraic closure \,$\clfld$. 

Let $U$ be an open subvariety of $X$. 
The geometric Galois coverings 
that are \emph{finite and \'etale over U} and defined over $\fld$ 
are classified by 
\begin{align*}
	\fundGO{U} :=& 
	\varprojlim_{\substack{V \ra U \\ \textrm{ geo.\ fin.\ \'et.\ Galois}/\fld}} \Gal\pn{K_V|K_U} 
\end{align*} 
and the \emph{geometric \'etale fundamental group} \,$\fundGgeo{U}$\, 
is the inverse limit where in addition the field of definition 
ranges over finite field extensions \,$\erwfld | \,\fld$\, in \,$\clfld$.

\section{Related Problems}
\label{sec:Open Problems}

The first problem is the central construction that is used 
in order to approach fundamental groups via duality. 

\bPbl[Relation between \'etale fundamental group and Tannakian fundamental group] 
\label{pbl: group fctr of rel Cartier div}
Let $X$ be a projective curve and \,$U$ an affine subvariety. 
The goal is to find a Tannakian category $\Cat$ 
whose Tannakian fundamental group 
is (as close as possible to) the geometric fundamental group of \,$U$. 
\ePbl 

Let us consider the abelian case first. 
The abelianized geometric fundamental group of $U$ 
is associated to the Cartier dual 
of the completion of the functor of principal divisors $\HDivsf{X}{S}$ 
on $X$ with support in $S = X \setminus U$: 
\begin{align} 
	R \,\lmt\, \HDivsfc{X}{S}(R) := \Gm\pn{\sO_X(U) \tens_k R} \big/ \Gm(R) 
	\laurin 
	\label{HDiv} \mytag{H}
\end{align} 
Here $R$ ranges over finite commutative $\fld$-rings. 
If we extend this assignment to non-commutative rings as arguments, 
the values are only in homogeneous spaces.  
As for fundamental groups, we will need to fix a base point $\bpt$ on $U$ 
in order to settle a group structure. 
Consider the elements of $\Gm\bigpn{\sO(U) \tens R}$ 
as morphisms from $U$ to the Weil restriction $\Lin_R := \WeilRes{R}{\fld} \Gm$. 
The base point $\bpt$ allows us to pick 
for every residue class \,$\Del \in \Gm\bigpn{\sO(U) \tens_k R} \big/ \Gm(R)$\, 
a distinguished representative \,$\del \in \Del$, 
namely the morphism \,$\del: U \lra \Lin_R$\, that maps \,$\bpt \in U$\, 
to \,$1 \in \Lin_R$. 
For $V$ a $\fld$-vector space, setting $R = \End(V)$ yields \,$\Lin_R = \Gl_V$, 
the group of automorphisms on $V$. 
In this way we replace quotient spaces by subgroups 
and obtain the functor  
\begin{align} 
	V \,\lmt\, \Autfs{X,\bpt}{\,S}(V) := \Mor\bigpn{\pn{U,\bpt}, \pn{\Gl_V, 1}}
	\label{Gbm Mor} \mytag{M} 
\end{align} 
which has values in maps to automorphisms of a vector space. 
We will refer to $\Autfs{X,\bpt}{\,S}$ as the 
\emph{functor of based morphisms from $U$ to automorphisms}.  
Similarly, we define 
\begin{align} 
	V \,\lmt\, \Autf_{X,\bpt}(V) := \RatMap\bigpn{\pn{X,\bpt}, \pn{\Gl_V, 1}}
	\label{Gbm ratlMap} \mytag{N} 
\end{align} 
which we will refer to as the 
\emph{functor of based rational maps from $X$ to automorphisms}.  
As we will show in Section \ref{sec:TannakaDuality} and Proposition \ref{Gbm fmlGroup}, 
the functors $\Autf_{X,\bpt}$ and $\Autfs{X,\bpt}{\,S}$ encode 
categories \,$\Cat\bigpn{\Autf_{X,\bpt}}$ resp.\ \,$\Cat\bigpn{\Autfs{X,\bpt}{\,S}}$ 
of representations of affine groups. 
We denote by $\Autfc_{X,\bpt}$ resp.\ $\Autfcs{X,\bpt}{\,S}$ 
the equivalence classes of those functors that define equivalent categories. 

Then $\Cat\bigpn{\Autfs{X,\bpt}{\,S}}$ will turn out to be 
the desired Tannakian category: 
the geometric fundamental group of $U$ is 
the pro-sub-completion (see Notations \ref{pro-sub-completion}) 
of the Tannaka fundamental group of $\Cat\bigpn{\Autfs{X,\bpt}{\,S}}$, 
cf.\ Problem \ref{genl FundGroup of Affine}.  


The next problem 
is a generalization of Theorem \ref{FundGroup-of-A1} 
(more precisely of Theorem \ref{FundGroup-of-Affine}), 
but this will become clear only in the course of these notes. 

\bPbl[Explicit description of absolute Galois groups and fundamental groups]
\label{genl FundGroup of Affine}
Let $U$ be an affine variety and $X$ a projective compactification, 
over a finite field $\fld$ with an algebraic closure $\clfld$. 
The goal is to describe the geometric fundamental group of \,$U$ and 
the geometric absolute Galois group of the function field $K_X$ 
in terms of $\Autfc_{X,\bpt}$ 
and $\Autfcs{X,\bpt}{\,S}$. 
The question is, in which generality the following descriptions hold:  
\begin{align*} 
\fundG\bigpn{U_{\clfld},\bpt} 
	&= \ProSub\Bigpn{\bigpn{\Autfcs{X,\bpt}{\,S}}^{\vee}\bigpn{\clfld} } \\ 
\Gal\bigpn{K_{\ol{X}}^{\sep} \big| K_{\ol{X}}} 
	&= \bigcup_{\bpt \in X(\clfld)}  \ProSub\Bigpn{\bigpn{\Autfc_{X,\bpt}}^{\vee}\bigpn{\clfld} } 
	\laurin 
\end{align*} 
\ePbl
\noindent 
Here $\Autfc_{X,\bpt}$ and $\Autfcs{X,\bpt}{\,S}$ 
are equivalence classes of functors of representations of affine groups, 
their Tannaka duals (denoted by $\lul^{\vee}$, see Theorem \ref{Tannaka duality}) 
are affine groups, 
which we can take $\clfld$-valued points of. 
Of the resulting group we consider the pro-sub-completion 
(see Notations \ref{pro-sub-completion}). 

\bPbl[Recover coordinate ring of a variety from its fundamental group] 
\label{Recover CoordRing from FundGroup}
Conversely, suppose 
we are given a profinite group $\Pi$. 
Then we are looking for an affine variety $U$ such that \,$\Pi$\, 
is the geometric fundamental group of \,$U$: 
\[ \Pi = \fundGgeo{U} 
:= \varprojlim_{\substack{V\,|\,U \\ \textrm{ geometric finite \'etale Galois}}} 
   \Gal\bigpn{V\,\big|\,U} 
\] 
where \,$\clfld$\, is an algebraic closure of the field of definition \,$\fld$\, of \,$U$. 
\ePbl 

If we find an affine group $\Gp$ with \,$\ProSub\bigpn{\Gp\pn{\clfld}} = \Pi$\, 
(i.e., the quotients of $\Pi$ 
are the subgroups of the quotients of $\Gp(\clfld)$), 
then the Tannaka dual \,$\Reprf = \Gpd$\, of $\Gp$ 
helps us to recover the coordinate ring \,$\sO(U)$\, of \,$U$: 
The problem is reduced to finding a finitely generated $\fld$-algebra $A$ 
with an augmentation \,$\aug: A \ra \fld$\, such that 
the following functor from $\fld$-vector spaces to sets 
\begin{align*}
	V \;&\lmt\; \ker\Bigpn{\Gl_V\pn{A} \xra{\Gl_V(\aug)} \Gl_V\pn{\fld}} 
\end{align*}
is equivalent to the functor of representations of $\Gp$ 
(see Definition \ref{representation functor} 
for the notion of a \emph{functor of representations}). 
This is the inversion of Problem \ref{genl FundGroup of Affine}, 
i.e.\ we use that the functor \,$\Autfs{X,\bpt}{S}: \Veck \lra \Set$\, 
from Definition \ref{group of based maps} (cf.\ (\ref{Gbm Mor})) 
is a representative for the Tannaka dual of the affine group 
that belongs to the geometric fundamental group of $\pn{X \setminus S, \bpt}$.

\section{Tannaka Duality}
\label{sec:TannakaDuality}

The key tool for our construction of universal affine groups 
is duality theory for affine groups. 
In the abelian case \emph{Cartier duality} was the appropriate theory. 
In the non-abelian case we want to replace this by \emph{Tannaka duality}. 
Unfortunately, the Tannakian theory involves a large pile of 
categorical concepts. 
Here we refer to \cite[II]{DelMil_Tannaka} as a main reference. 
In this section we aim to give a reformulation of Tannaka duality 
that is convenient for our purpose, 
e.g.\ formally analogous to Cartier duality would be helpful. 
Especially, what we need is an \emph{anti-equivalence} of affine groups 
to some other category that is usable and handy 
(in particular \emph{not} a ``category of categories''). 
Let $\fld$ be a field. 

\bDef
\label{GL_V}
For a $\fld$-vector space $V$ we let $\Gl_V$ be the $\fld$-group functor 
that assigns to a $\fld$-algebra $R$ the group of $R$-linear automorphisms 
of $V \tens R$: 
\[ \Gl_V: R \lmt \Aut\pn{V \tens R} 
\laurin 
\]
\eDef 

\bDef 
\label{Rep_G}
Let $\Gp$ be an affine group. 
Then \,$\Rep_{\Gp}$ denotes the \emph{category of representations of $\Gp$}: 
the objects are $G$-modules $(V,\rho)$, i.e.\ vector spaces $V$ 
that admit a representation \,$\rho \in \Homgk\pn{\Gp, \Gl_V}$, 
and morphisms between $\Gp$-modules $(V,\rho_V)$ and $(W,\rho_W)$ 
are $k$-linear maps $\phi: V \ra W$ 
that are $\Gp$-equivariant, 
i.e.\ satisfying the condition 
\,$\phi \circ \rho_V(g) = \rho_W(g) \circ \phi$\, for all $g \in \Gp$. 
\eDef 

\bDef 
\label{Repf_G}
For an affine group $\Gp$ 
we let \,$\Repf_G: \Veck \lra \Set$\, 
be the following functor 
on finite dimensional $\fld$-vector spaces with values in sets:  
\[ \Repf_{\Gp}: V \lmt \Homgk\pn{\Gp, \Gl_V} 
\laurin 
\]
\eDef 

\bDef 
\label{S-module}
Let $X$ be a $\fld$-functor. 
An $X$-module is a pair $(V,\alp)$, i.e.\ a $\fld$-vector space $V$ 
together with a map \,$\alp \in \Mor_{\Fctrk}\pn{X,\Gl_V}$, 
inducing the structure of $\Gpa(\alp)$-module on $V$, 
where $\Gpa(\alp) \subset \Gl_V$ is the subgroup generated by $\alp: X \ra \Gl_V$. 
A homomorphism of $X$-modules from $(V,\alp_V)$ to $(W,\alp_W)$ 
is a $k$-linear map $\phi: V \ra W$ 
that is $X$-equivariant, i.e.\ satisfying 
\,$\phi \circ \alp_V(x) = \alp_W(x) \circ \phi$\, for all $x \in X$. 
The category of $X$-modules is denoted by $X$-$\Mod$. 
\eDef 

\bDef 
\label{action functor}
Let $X$ be a $\fld$-functor. 
A \emph{functor of actions generated by $X$} (or \emph{of $X$-actions}) 
is a functor \,$\ul{A}: \Veck \lra \Set$\, 
that assigns to a finite dimensional $\fld$-vector space $V$ 
a subset \,$\ul{A}(V) \subset \Mor_{\Fctrk}\pn{X,\Gl_V}$.  
A homomorphism of functors of actions from $\ul{A}$ with generator $X$ 
to $\ul{B}$ with generator $Y$ 
is by definition a natural transformation \,$\ul{A} \ra \ul{B}$\, 
that is induced by a map $\morphism \in \Mor_{\Fctrk}(Y, X)$. 
Associated to a functor of actions $\ul{A}$ 
is a \emph{category of actions} $\Cat(\ul{A})$: 
if $X$ is the generator of actions, 
this is the full subcategory of $X$-$\Mod$ 
whose objects are $X$-modules $(V,\alp)$, where $V$ is a finite dimensional 
$\fld$-vector space and \,$\alp \in \ul{A}(V)$. 
\eDef 

\bRmk 
\label{ReprFctr ReprCat}
(a) Let $\Gp$ be an affine group. 
Then $\Repf_{\Gp}$ is a functor of $\Gp$-actions, 
and we have \,$\Cat\bigpn{\Repf_{\Gp}} = \Rep_{\Gp}$. 

(b) Let $X$ and $Y$ be $\fld$-functors, 
$\ul{A}$ a functor of $X$-actions and 
$\ul{B}$ a functor of $Y$-actions. 
A map \,$\morphism: \Mor_{\Fctrk}(Y, X)$\, 
induces a homomorphism of functors of actions \,$\ul{A} \ra \ul{B}$\, 
only if \,$\morphism^*: \Mor_{\Fctrk}(X,\Gl_V) \lra \Mor_{\Fctrk}(Y,\Gl_V)$, 
$\alp \lmt \alp \circ \morphism$\, 
maps \,$\ul{A}(V) \ra \ul{B}(V)$\, for every finite dimensional vector space $V$. 
\eRmk 

\bLem 
\label{Hom of Repf}
If $G$ and $H$ are affine groups, 
then a map \,$\morphism \in \Mor_{\Fctrk}(H, G)$\, 
induces a homomorphism \,$\Repf_{\morphism}: \Repf_{G} \lra \Repf_{H}$\, 
if and only if \,$\morphism$\, is a homomorphism of \,$\fld$-groups, 
i.e.\ \,$\morphism \in \Homgk(H, G)$. 
\eLem 

\bPf 
If \,$\morphism \in \Homgk(H, G)$, then obviously 
\,$\morphism^*: \alp \lmt \alp \circ \morphism$\, 
maps \,$\Homgk(G,\Gl_V) \lra \Homgk(H,\Gl_V)$. 

The converse direction can be seen as follows: 
Suppose first that $G$ is an affine algebraic group. 
Then $G$ has a faithful representation $\rho \in \Homgk(G,\Gl_V)$,  
and $\rho$ has an inverse \,$\sig \in \Homgk(\im\rho, G)$\, 
with \,$\sig \circ \rho = \id_G$. 
Hence \,$\rho \circ \morphism \in \Homgk(H,\Gl_V)$\, 
if and only if \,$\sig \circ \rho \circ \morphism = \morphism \in \Homgk(H, G)$. 

In general, every affine group $G$ is pro-algebraic 
(see \cite[II, Cor.\ 2.7]{DelMil_Tannaka}). 
This means that $G = \varprojlim G_i$ is an inverse limit of affine algebraic groups $G_i$, 
and we can apply the above argument to every quotient $G_i$. 
\ePf

\bigskip 
The following result is the main theorem for Tannaka duality. 
For notation and proof we refer to 
\cite[Theorem 2.11]{DelMil_Tannaka}: 

\bThm[Deligne, Milne]
\label{Deligne Milne}
Let $(\Cat,\tens)$ be a rigid abelian tensor category such that \,$\fld = \End(\unit)$, 
and let \,$\oma: \Cat \ra \Veck$\, be an exact faithful $\fld$-linear tensor functor. 
Then, 

\noindent 
\begin{tabular}{rl}
	{\rm (a)} & the functor $\Autf^{\tens}(\oma)$ of \,$\fld$-algebras is represented by an affine group $\Gp$\laurink \\ 
	{\rm (b)} & the functor $\Cat \ra \Rep_{\Gp}$ defined by $\oma$ is an equivalence 
	of tensor categories\laurin
\end{tabular}
\eThm 

\bDef 
\label{representation functor} 
A \emph{functor of representations} is by definition a functor of actions \,$\ul{R}$\, 
for which the associated category $\Cat(\ul{R})$ and the natural forgetful functor 
$\oma: \Cat(\ul{R}) \lra \Veck$\, satisfy the conditions of Theorem \ref{Deligne Milne}. 

Two functors of actions $\ul{R}$ and $\ul{S}$ are said to be \emph{equivalent}, 
if there exists an equivalence of associated categories \,$\Cat(\ul{R}) \iso \Cat(\ul{S})$. 

A \emph{representation functor} $\Reprf$ is by definition an equivalence class $[\ul{R}]$ 
of a functor of representations. 
Homomorphisms of representation functors are those induced by 
homomorphisms of functors of representations by passing to equivalence classes. 
(This definition is justified by Proposition \ref{hom of repr functors}) below. 
\eDef 

\bRmk
\label{functor of repr criterion}
A functor of actions $\ul{A}$ 
is a functor of representations 
if and only if it is equivalent to a functor of representations $\Repf_{\Gp}$ 
of some affine group $\Gp$. 
(Immediate from Definition \ref{representation functor} and Theorem \ref{Deligne Milne}.)
\eRmk 

\bPrp[Homomorphisms of Representation Functors] 
\label{hom of repr functors}
A homomorphism of functors of representations \,$h: \ul{R} \ra \ul{S}$\, 
induces a well defined map on equivalence classes \,$[h]: [\ul{R}] \ra [\ul{S}]$. 
\ePrp 

\bPf 
Let $\ul{R}$ be a functor of representations with generator $X$. 
Then the associated category $\Cat\pn{\ul{R}}$ is exactly determined 
by the various subgroups $\Gpa(\alp)$ of $\Gl_V$ generated by 
\,$\alp \in \ul{R}(V) \subset \Mor_{\Fctrk}(X,\Gl_V)$, 
where $V$ runs through finite dimensional $\fld$-vector spaces, 
cf.\ Remark \ref{G(X) for reprf}. 
Therefore passing to equivalence classes \,$\ul{R} \lmt [\ul{R}]$\, 
means to replace the $\alp$ by $\Gpa(\alp)$. 

Let $\ul{S}$ be another functor of representations with generator $Y$, 
and \,$h: \ul{R} \ra \ul{S}$\, a homomorphism. 
Then $h$ is given by $\morphism^*: \alp \lmt \alp \circ \morphism$ 
for some morphism $\morphism: Y \ra X$. 
Passing to equivalence classes $h \lmt [h]$ means to replace 
the assignments \,$\alp \lmt \alp \circ \morphism$\, 
by \,$\Gpa(\alp) \lmt \Gpa(\alp \circ \morphism)$, 
which does not depend on the choices of $\alp$ or $\morphism$, 
but only on the subgroups they generate and the association between them. 
This gives a well defined map \,$[\ul{R}] \ra [\ul{S}]$. 
\ePf 

\bRmk
\label{G(X) for reprf} 
The main Theorem \ref{Deligne Milne} says 
that every representation functor $[\ul{R}]$ 
has a representative $\Repf_{G}$ for some affine group $G$. 
In particular this implies that the data $\bigst{\Gpa(\alp) \,\big|\, \alp \in \ul{R}(V)}_V$, 
which determines the category of representations $\Cat(\ul{R})$, 
is equivalent to the data $\bigst{\Homgk(G,\Gl_V)}_V$ for some affine group $G$.  
In other words, the $\Gpa(\alp)$ are given by the homomorphic images of some universal group $G$. 
(By definition we assume here that $\Cat(\ul{R})$ 
satisfies the conditions of Theorem \ref{Deligne Milne}.) 
\eRmk 

These constructions allow us to reformulate Tannaka duality 
in the following way: 

\bThm[Tannaka Duality] 
\label{Tannaka duality}
The following assignment establishes an anti-equivalence of categories 
between affine groups and representation functors: 
\begin{align*}
	\bigpn{\textrm{affine $\fld$-groups}} & \iso \bigpn{\textrm{representation functors}} \\
	\Gp & \lmt \bigbt{\Repf_{\Gp}} 	 
\end{align*}
The converse direction is given for a representation functor $\ul{R}$ by 
\[ \bigbt{\ul{R}} \lmt \Autf^{\tens}(\oma) 
\] 
where \,$\oma: \Cat(\ul{R}) \lra \Veck$\, is the natural forgetful functor. 
\eThm 

\bPf 
The association between objects follows from Theorem \ref{Deligne Milne}. 
For two representation functors $\Reprf$ and $\Reprff$ 
we can therefore choose representatives $\Repf_{G}$ and $\Repf_{H}$. 
The homomorphisms of representation functors from $\Reprf$ to $\Reprff$ 
are then induced by 
the homomorphisms $\Hom_{\FctrAct}\pn{\Repf_{G},\Repf_{H}}$ 
of functors of representations from $\Repf_{G}$ to $\Repf_{H}$, 
which correspond exactly to the homomorphisms $\Homgk\pn{H,G}$ 
of affine groups from $H$ to $G$, 
due to Lemma \ref{Hom of Repf}. 
\ePf

\section{Recovering a Group from its Representations}
\label{sec:Recover from Repr}

Given a representation functor $\Reprf$ over a field $k$, 
the goal is to describe the affine group 
that is given by the Tannaka dual $\Gp := \Reprfd$ 
in terms of its coordinate ring and group structure. 
We will apply the method of recovering a group from its representations, 
as described in \cite[No.~2.5]{Springer}, 
using Tannaka's Theorem \cite[2.5.3]{Springer} as a main ingredient. 

\bPnt 
\label{Springer Recovery}
We sketch the method from \cite[No.~2.5]{Springer} in our own words: 
Let $\Gp$ be an (unknown) reduced affine group. 
Then we can recover its algebra $\sO(\Gp)$ and its coalgebra structure 
$\comult: \sO(\Gp) \lra \sO(\Gp) \tens \sO(\Gp)$ 
from the category of finite dimensional $\Gp$-modules. 
Here objects are vector spaces $V$ that admit a representation $\rho: \Gp \ra \Aut(V)$, 
and morphisms between $\Gp$-modules $V$ and $W$ are $k$-linear maps $\phi: V \ra W$ 
that are $\Gp$-equivariant, 
i.e.\ satisfying the condition 
\,$\phi \circ \rho_V(g) = \rho_W(g) \circ \phi$\, for all $g \in \Gp$. 
Consider the $k$-algebra 
\[ \sA = \varinjlim_{(V,\,\rho)} \,\End(V)^{\vee}
\]
the direct limit of endomorphisms of vector spaces $V$, 
where $(V,\rho)$ ranges over all finite dimensional $\Gp$-modules. 
This system is directed by pairs $\pn{\phi,\psi}$ of homomorphisms of $\Gp$-modules 
\,$\phi: V \ra W$\, and \,$\psi: W \ra V$\, 
that induce inverse isomorphisms \,$V/\ker(\phi) \iso W/\ker(\psi)$\, 
for which the $\Gp$-action on $\ker(\phi)$ is trivial. 
The transition map \,$\End(V)^{\vee} \lra \End(W)^{\vee}$ is defined by 
\,$v \tens f \lmt \phi(v) \tens \psi^{\vee}(f)$\, 
for \,$v \in V$\, and \,$f \in V^{\vee}$. 
For \,$h := \psi^{\vee}(f)$ this means \,$v \tens \phi^{\vee}(h) \lmt \phi(v) \tens h$.  
Let \,$a_V: \End(V)^{\vee} \lra \sA$\, be the canonical maps. 
The algebra multiplication of $\sA$ 
for $\alp \in \End(V)^{\vee}$ and $\bet \in \End(W)^{\vee}$ 
is 
\[ a_V(\alp) \cdot a_W(\bet) := a_{V \oplus W} \bigpn{\pr_1^*\alp + \pr_2^*\bet} 
\in a_{V \oplus W} \bigpn{\End\pn{V \oplus W}^{\vee}} 
\laurin 
\] 
The coalgebra structure on $\comult_{\sA}: \sA \lra \sA \tens \sA$ is defined as follows: 
For a \,$\Gp$-module $V$ choose a basis $\st{e_1, \ldots, e_n}$, 
and let $\st{e_1^{\vee}, \ldots, e_n^{\vee}}$ be the corresponding cobasis. 
Then \,$\id_{V} = \sum_i e_i^{\vee} \tens e_i \in V^{\vee} \tens V = \End(V)$ 
is the identity operator on $V$. 
For $v \in V$ and $f \in V^{\vee}$ we set 
\,$\comult_{V}\pn{v \tens f} = v \tens \id_{V} \tens f 
= \sum_i \pn{v \tens e_i^{\vee}} \tens \pn{e_i \tens f} 
\in \End(V)^{\vee} \tens \End(V)^{\vee}$, 
and we define \,$\comult_{\sA}\bigpn{a_V\pn{v \tens f}} 
= \pn{a_V \tens a_V} \bigpn{\comult_{V}\pn{v \tens f}}$.  
With \,$E_{ij} = a_V\pn{e_i \tens e_j^{\vee}}$ this can be written as 
\[ \comult_{\sA}: E_{ij} \lmt \sum_k E_{ik} \tens E_{kj} 
\laurin 
\]
According to \cite[Thm.\ 2.5.7]{Springer} there is an isomorphism of $k$-algebras 
\begin{align} 
	\sA/\Nil(\sA) &\iso \sO(\Gp)  \label{algebra of reps} \mytag{A}
\end{align} 
where $\Nil(\sA)$ is the nil-radical of $\sA$. 
Moreover, $\comult_{\sO(\Gp)}$ is induced by $\comult_{\sA}$. 
\ePnt 

\bRmk[Deduction of $\sO(\Gp)$] 
\label{deduktion of O(G)}
Let $\Gp$ be an affine group, 
$\Rep_{\Gp}$ its category of representations 
and $\oma: \Rep_{\Gp} \lra \Veck$ the forgetful functor. 
Let $\End(\oma)$ be the natural transformations of $\oma$ on itself. 
An element of $\End(\oma)$ is a family $\st{\Tha_V}_{V \in \Gp\textrm{-}\Mod}$ 
of compatible endomorphisms \,$\Tha_V \in \End(V)$, 
i.e.\ satisfying \,$\phi \circ \Tha_V = \Tha_W \circ \phi$\, 
for every homomorphism of $\Gp$-modules \,$\phi: V \ra W$. 
If $\Gp$ is finite, 
the affine algebra of $\Gp$ is the $k$-linear dual to $\End(\oma)$, 
see \cite[9.e, 9.42]{Milne_algGroups}: 
\[ \sO(\Gp) = \End(\oma)^{\vee}
\laurin 
\]
In order to write $\End(\oma)$ as an inverse limit, 
we need for every homomorphism \,$\phi: V \ra W$ 
a reciprocal homomorphism \,$\psi: W \ra V$ 
that induces an inverse to \,$[\phi]: V / \ker(\phi) \iso \im(\phi)$. 
This is equivalent to the choice of splittings \,$V \cong \ker(\phi) \dsum \im(\phi)$\, 
and \,$W \cong \im(\phi) \dsum \coker(\phi)$. 
Such splittings are not canonical, and instead of choosing one 
we take all of them by considering an inverse system 
with pairs $\pn{\phi,\psi}$ as transition functions from $V$ to $W$ 
as in \Point \ref{Springer Recovery}. Then 
\begin{align*}
	\End(\oma) &= \varprojlim_{(V,\,\rho)} \End(V) \\ 
	\sO(\Gp) &= \varinjlim_{(V,\,\rho)} \End(V)^{\vee} 
	\laurin 
\end{align*}
\eRmk

\bPnt 
\label{Repr from CartierDual}
Let $\Reprf$ be a representation functor, 
such that its Tannaka dual $\Gp := \Reprfd$ is a unipotent affine group. 
We aim to describe $\sO(\Gp)$. 
By Lemma \ref{limit computation} 
the algebra \,$\sA = \varinjlim \End(V)^{\vee}$\, 
from \ref{Springer Recovery} is  
\,$\sA = \varinjlim \Homgk\bigpn{\Gp, \ul{\End}(V)^{\vee}}$. 
Identifying faithful representations on vector spaces of the same dimension 
and denoting \,$\gl_n = \ul{\Lie}\pn{\Gl_n} \cong \ul{\End}\pn{k^{n}}^{\vee}$  we get 
\,$\sA = \varinjlim_n \;\Sym\bigpn{\Homgk\pn{\Gp,\gl_n}}$. 
Here the algebra multiplication of $\sA$ from \ref{Springer Recovery} 
is expressed in terms of the $\Sym$-functor, 
in order to consider every summand \,$\Sym\bigpn{\Homgk\pn{\Gp,\gl_n}}$\, 
in the direct limit as an algebra on its own. 
By \Point \ref{Springer Recovery} equation(\ref{algebra of reps}) 
the algebra of $\Gp$ is \,$\sO(\Gp) = \sA / \Nil(\sA)$. 
Taking into account that $\sA$ does not have any nilpotent elements in this situation, 
we obtain 
\begin{align} 
	\sO(\Gp) 
	=\, \varinjlim_n \;\Sym\bigpn{\Homgk\pn{\Gp,\gl_n}} 
	\label{O(L)} \mytag{O} 
	\laurin 
\end{align}
\ePnt 

\bLem 
\label{limit computation}
Let $\Gp$ be a unipotent affine group. 
There exists an isomorphism of \,$k$-algebras 
\begin{align*} 
	\varinjlim_{(V,\,\rho)} \,\End(V)^{\vee} 
	\;&\iso\; \varinjlim_{(V,\,\rho)} \;\Homgk\bigpn{\Gp, \ul{\End}(V)^{\vee}} 
\end{align*} 
where $(V,\rho)$ ranges over all finite dimensional $\Gp$-modules. 
\eLem 

\bPf 
We will show the assertion by finding monomorphisms of $\Gp$-modules, 
compatible with the algebra structures, in either direction. 

{\bf Step 1:} \;$\sA \linj \Hom(\Gp,\sA)$ \\ 
Since $\Gp$ is unipotent, it admits a non-trivial homomorphism to $\Ga$ 
(see \cite[IV, \S2, 2.1]{DG}). 
Such a \,$\pi \in \Homgk\pn{\Gp, \Ga}$ 
induces a monomorphism of $\Gp$-modules 
\begin{align}
	V \linj \Homgk\pn{\Gp,\ul{V}} \laurink \hspace{35pt} v \lmt \emb_v \circ \pi 
	\label{mono} \mytag{U} 
\end{align}
where \,$\emb_v \in \Homgk\pn{\Ga,\ul{V}}$ denotes the homomorphism 
sending \,$1 \lmt v$\, for $v \in V$, 
and the $\Gp$-module structure of $\Homgk\pn{\Gp,\ul{V}}$ is inherited from $V$. 
Replacing $V$ by $\End(V)^{\vee}$ 
and taking the limit over finite dimensional $\Gp$-modules 
gives the first monomorphism 
\begin{align*} 
\varinjlim_{V} \,\End(V)^{\vee} 
&\linj\; \varinjlim_{V} \;\Homgk\bigpn{\Gp, \ul{\End}(V)^{\vee}} 
\laurin 
\end{align*}

{\bf Step 2:} \;$\Hom(\Gp,\sA) \linj \sA$ \\ 
The existence of a non-trivial \,$\pi \in \Homgk\pn{\Gp, \Ga}$ 
and the surjectivity of the assignment 
\,$\Homgk\pn{\Ga,\ul{V}} \lra V$, \,$\emb_v \lmt v$\, 
ensure that we have an epimorphism of $k$-vector spaces 
\;$ \chi: \Homgk\pn{\Gp,\ul{V}} \lsur V 
$. 
Left-exactness of the $\Hom$-Functor yields a monomorphism 
\[ \Homk\pn{V, H} \linj \Homk\pn{\Homgk\pn{\Gp,\ul{V}}, H} 
\laurink \hspace{35pt} f \lmt f \circ \chi 
\laurin 
\]
Setting \,$H := \Homgk\pn{\Gp,\ul{V}}$\, and 
composing with a preceding reordering of arguments 
we obtain a monomorphism 
\begin{align*} 
	\Homgk\bigpn{\Gp, \ul{\Hom}_k(V, V)} 
	\hspace{3pt} \iso & \;\Homk\bigpn{V, \Homgk\pn{\Gp,\ul{V}}} \\ 
	\linj & \;\Homk\bigpn{\Homgk\pn{\Gp,\ul{V}}, \Homgk\pn{\Gp,\ul{V}}} 
	\laurin 
\end{align*}
We have thus achieved a monomorphism 
\begin{align*} 
	\Homgk\bigpn{\Gp, \ul{\End}(V)} 
	\linj & \;\End\bigpn{\Homgk\pn{\Gp,\ul{V}}} 
	\laurin 
\end{align*}
Since every affine group is pro-algebraic (see \cite[II, Cor.\ 2.7]{DelMil_Tannaka}), 
we can write \,$\Gp = \varprojlim_i \Gp_i$, where the $\Gp_i$ are algebraic groups.  
Then 
\begin{align}
	\varinjlim_V \;\Homgk\bigpn{\Gp, \ul{\End}(V)}  
	\hspace{3pt}=& \;\varinjlim_{i,V} \;\Homgk\bigpn{\Gp_i, \ul{\End}(V)} 
	\notag \\ 
	\linj & \;\varinjlim_{i,V} \;\End\bigpn{\Homgk\pn{\Gp_i,\ul{V}}} 
	\label{Hom limit} \mytag{V} 
	\laurin  
\end{align}
Set \,$H_{i,V} := \Homgk\pn{\Gp_i,\ul{V}}$, these are finite dimensional $\Gp$-modules. 
Then $\st{H_{i,V}}_{i,V}$ is a directed subset of $\st{V}_{(V,\,\rho)}$, 
and the subsystem $\st{\End\pn{H_{i,V}}}_{i,V}$ 
satisfies the Embedding Condition $\lrpn{\bigodot}$ 
from Definition \ref{Embedding conditioon} 
in $\st{\End\pn{V}}_{(V,\,\rho)}$: 
Suppose \,$H_{i,X} < V$, \,$H_{j,Y} < V$, 
\;$\xi \in \End\pn{H_{i,X}}$, \,$\eta \in \End\pn{H_{j,Y}}$, 
and \,$\trans_{V,\,H_{i,X}}\pn{\xi} = \trans_{V,\,H_{j,Y}}\pn{\eta}$. 
As $\st{\Gp_i}_i$ is a directed inverse system,  
there is an index \,$k$ such that 
we have maps \,$\Gp_k \ra \Gp_i$, \,$\Gp_k \ra \Gp_j$. 
Moreover there is a monomorphism \,$V \inj H_{k,V}$\, by Formula (\ref{mono}).  
Then \,$H_{i,X} < H_{k,V}$, \,$H_{j,Y} < H_{k,V}$, 
and \,$\trans_{H_{k,V},H_{i,X}}\pn{\xi} 
= \pn{\trans_{H_{k,V}, V} \circ \trans_{V,\,H_{i,X}}}\pn{\xi}
= \pn{\trans_{H_{k,V}, V} \circ \trans_{V,\,H_{j,Y}}}\pn{\eta}
= \trans_{H_{k,V},\,H_{j,Y}}\pn{\eta}$. 
By Lemma \ref{Embedding lemma} we get a monomorphism 
\begin{align}
	\varinjlim_{i,V} \,\End\pn{H_{i,V}} \;\linj\; \varinjlim_{V} \,\End\pn{V}
	\label{End limit} \mytag{W}
	\laurin 
\end{align} 
Composition of the monomorphisms (\ref{Hom limit}) and (\ref{End limit}) 
and the canonical isomorphism 
\[ \End\pn{V} \cong V^{\vee} \tens V \cong V \tens V^{\vee} 
\cong \pn{V^{\vee} \tens V}^{\vee} \cong \End(V)^{\vee} 
\] 
yield the desired monomorphism 
\[ \varinjlim_V \;\Homgk\bigpn{\Gp, \ul{\End}(V)^{\vee}} \;\linj\; 
   \varinjlim_{V} \;\End(V)^{\vee} 
\laurin 
\]
\hfill 
\ePf 

\bigskip
The following definition and lemma are taken from 
\cite[Definition 3.1 and Lemma 3.2]{Rubin}. 

\bDef
\label{Embedding conditioon} 
Let $\bfX = \pn{X_{\alp}, \trans_{\alp \bet}, \pn{A,\leq}}$ a direct system 
in a category $\Cat$ that admits direct limits, 
and $N \subset A$ a directed subset. 
Set $\bfY = \pn{X_{\nu}, \trans_{\nu \mu}, \pn{N,\leq}}$. 
The \emph{embedding condition} of \,$\bfY$\, in \,$\bfX$\, is as follows: 

\bigskip 
\noindent 
\begin{tabular}{rl}
	$\lrpn{\bigodot}$ & If \,$\lam, \mu \in N$, $\alp \in A$, 
	$\lam \leq \alp$, $\mu \leq \alp$, 
	\;$y \in X_{\lam}$, \,$x \in X_{\mu}$, \\
	& and \,$\trans_{\alp \lam}\pn{y} = \trans_{\alp \mu}\pn{x} \in X_{\alp}$, \\ 
	& then there exists \,$\nu \in N$\, 
	such that \,$\lam \leq \nu$, \,$\mu \leq \nu$, \\ 
	& and \,$\trans_{\nu \lam}\pn{y} = \trans_{\nu \mu}\pn{x} \in X_{\nu}$. 
\end{tabular} 
\eDef 

\smallskip 
\bLem 
\label{Embedding lemma} 
With the notation as in Definition \ref{Embedding conditioon}, 
if \,$\bfY$ satisfies the embedding condition $\lrpn{\bigodot}$ in $\bfX$, 
then the induced map on direct limits \,$\varinjlim \bfY \lra \varinjlim \bfX$\, 
is a monomorphism:
\[ \varinjlim_{\nu \in N} X_{\nu} \linj\, \varinjlim_{\alp \in A} X_{\alp}
\laurin 
\]
\eLem 

\bPf 
\cite[Lemma 3.2]{Rubin}. 
\ePf

\section{Functor of Based Rational Maps to Linear Groups}
\label{sec:CategoryAut}

Let $\pn{X,\bpt}$ be a projective variety 
with a base point over a field $k$. 
 
\bDef 
\label{group of based maps}
Let \,$\AutfXb$\, be the functor that assigns to 
a $\fld$-vector space $V$ the set 
\begin{align*} 
\Autf_{X,\bpt}(V) 
&= \RatMap\bigpn{\pn{X,\bpt}, \pn{\Gl_{V}, 1}} \\
&= \ker\Bigpn{\Gl_V\pn{\sO_{X,\bpt}} \lra 
	\Gl_V\bigpn{\resfld\pn{\bpt}}} 
\end{align*} 
where \,$\sO_{X,\bpt}$\, is the local ring 
and \,$\resfld(\bpt)$\, the residue field at \,$\bpt$. 

Let \,$U \ni \bpt$\, be an open subvariety of $X$ containing $\bpt$ 
and set \,$S = X \setminus U$, 
then \,$\Autfs{X,\bpt}{\,S}$ is the functor defined by 
\begin{align*} 
	\Autfs{X,\bpt}{\,S}(V) 
	&= \Mor\bigpn{\pn{U,\bpt}, \pn{\Gl_{V}, 1}} \\
	&= \ker\Bigpn{\Gl_V\pn{\sO_{X}(U)} \lra 
		\Gl_V\bigpn{\resfld\pn{\bpt}}} 
	\laurin 
\end{align*} 
\eDef 

\bPrp 
\label{Gbm fmlGroup}
The functors \,$\Autf_{X,\bpt}$\, and \,$\Autfs{X,\bpt}{\,S}$\, 
are functors of representations of affine groups. 
\ePrp 

\bPf 
The associated categories $\Cat\pn{\Autf_{X,\bpt}}$ and $\Cat\pn{\Autfs{X,\bpt}{\,S}}$ 
satisfy the conditions (a)--(d) of \cite[Thm.\ 9.24]{Milne_algGroups}. 
Then the assertion follows by \,{\it loc.\ cit.} 
and Remark \ref{functor of repr criterion}. 
\ePf 

\bRmk 
\label{gbm iso}
An automorphism \,$\aut: X \iso X$ on $X$ induces 
isomorphisms of representation functors via pull-back \,$f \lmt \alp^*f$: 
\begin{align*}
	\alp^*: \Autf_{X,\bpt} & \iso \Autf_{X,\alp^{-1}\pn{\bpt}} \laurink & 
	\alp^*: \Autfs{X,\bpt}{\,S} & \iso \Autfs{X,\alp^{-1}\pn{\bpt}}{\,\alp^{-1}\pn{S}} 
	\laurin 
\end{align*}
\eRmk 

\bExm
\label{HDivf(gl) for A^1}
Consider \,$X = \Prj^1$, \,$S = \st{\infty}$\, 
and \,$U = \Prj^1 \setminus \st{\infty} = \Afn^1$.  
Then $\sO(U) = k[t]$, where $t$ is a local parameter of \,$0 \in \Prj^1$, and 
\begin{align*}
	\Autfs{\Prj^1, 0}{\st{\infty}}\bigpn{V} 
	&= \lrst{M = I_n + \sum_{i=1}^{r} A_i t^i \;\left|\; 
		\begin{array}{l}
			A_i \in \End(V) \\
			\det(M) \in \Gm\bigpn{k[t]}
		\end{array}
		\right.} 
		\laurin 
\end{align*}
\eExm 

\bExm 
\label{HDivf(gl) for Gm}
Let \,$X = \Prj^1$, \,$S = \st{0,\infty}$ 
and \,$U = \Prj^1 \setminus \st{0,\infty} = \Afn^1 \setminus \st{0}$. 
Then \,$\sO(U) = k[t, t^{-1}] = \bigcup_{\nu \in \Zint} t^{\nu} \,k[t]$, 
where $t$ is a local parameter of $0 \in \Prj^1$, and 
\begin{align*}
	\Autfs{\,\Prj^1, 1}{\st{0,\infty}}\bigpn{V} 
	&= \Biggst{\underbrace{t^{\nu}\Bigpn{I_n + \sum_{i=1}^{r} A_i \pn{t - 1}^i}}_M \;\Bigg|\; 
		\begin{array}{l}
			A_i \in \End(V), \,\nu \in \Zint \\
			\det(M) \in \Gm\bigpn{k[t,t^{-1}]}
		\end{array}
		} \\ 
	&\cong\, \Zint \tms \Autfs{\Prj^1, 0}{\st{\infty}}\bigpn{V} 
	\laurin 
\end{align*}
\eExm

\bPnt 
\label{Repr from CartierDual_Aut}
Now consider the representation functor \,$\Autfcs{X,\bpt}{\,S} = \bigbt{\Autfs{X,\bpt}{\,S}}$, 
and let \,$\LalbGs{X,\bpt}{\,S} := \bigpn{\Autfcs{X,\bpt}{\,S}}^{\vee}$ 
denote its Tannaka dual. 
Using Lemma \ref{limit computation_Aut} below 
we obtain similarly to \Point \ref{Repr from CartierDual} Equation \ref{O(L)}: 
\begin{align} 
	\sO\bigpn{\LalbGs{X,\bpt}{\,S}} 
	&=\, \varinjlim_n \;\Sym\Bigpn{\Mor\bigpn{\pn{U,\bpt},\pn{\gl_n, 0}}} 
	\label{O(L)_Aut} \mytag{O'} \\ 
	&=:\, \varinjlim_n \;\Sym\Bigpn{\Morpt\bigpn{U, \gl_n}} \notag 	
	\laurin 
\end{align}
\ePnt 

\bLem 
\label{limit computation_Aut}
Let $\pn{U,\bpt}$ be an affine $k$-scheme of finite type 
with a smooth base point. 
There exists an isomorphism of \,$k$-algebras 
\begin{align*} 
	\varinjlim_{(V,\,\alp)} \,\End(V)^{\vee} 
	\;&\iso\; \varinjlim_{(V,\,\alp)} 
	\;\Mor_{\Schk}\Bigpn{\bigpn{U,\bpt}, \bigpn{\ul{\End}(V)^{\vee},0}} 
\end{align*} 
where $(V,\alp)$ ranges over all finite dimensional $\pn{U,\bpt}$-modules. 
\eLem 

\bPf 
The proof is analogous to the proof of Lemma \ref{limit computation}. 
One has to replace \,$\Homgk\pn{\Gp,\ul{V}}$\, 
by \,$\Mor_{\Schk}\bigpn{\pn{U,\bpt}, \pn{\ul{V},0}}$, and 
a non-trivial homomorphism \,$\pi \in \Homgk\pn{\Gp, \Ga}$ 
by a non-trivial function 
\,$f \in \fm_{U,\bpt} \cong \Mor_{\Schk}\bigpn{\pn{U,\bpt}, \pn{\Afn^1,0}}$, 
where $\fm_{U,\bpt}$ is the maximal ideal of $\sO(U)$ at $\bpt$.  
Completion of local rings 
$\bigpn{\sO_{U,\bpt}, \fm_{U,\bpt}} \lra \bigpn{k[[t_1,\ldots,t_n]], (t_1,\ldots,t_n)}$ 
and the degree map \,$\deg: k[[t_1,\ldots,t_n]] \lra \Nat$\, 
allow us write the maximal ideal as a direct limit 
of finite dimensional $k$-vector spaces: 
\,$\fm_{U,\bpt} = \varinjlim_d \pn{\fm_{U,\bpt}}_{\leq d}$. 
Then 
\begin{align*}
	\Mor_{\Schk}\bigpn{\pn{U,\bpt}, \pn{\ul{V},0}} 
	&= \Homka\Bigpn{\bigpn{\Sym(V^{\vee}), (V^{\vee})}, \bigpn{\sO(U), \fm_{U,\bpt}}} \\ 
	&= \Homk\Bigpn{V^{\vee}, \varinjlim_d \,\pn{\fm_{U,\bpt}}_{\leq d}} \\ 
	&= \varinjlim_d \,\Homk\Bigpn{V^{\vee}, \pn{\fm_{U,\bpt}}_{\leq d}} 
\end{align*} 
is a direct limit of finite dimensional $k$-vector spaces as well. 
With this one can apply the same argumentation 
as in the proof of Lemma \ref{limit computation}. 
\ePf 

\bigskip 
\bExm 
\label{L from HDiv for A^1}
Consider \,$X = \Prj^1$, \,$S = \st{\infty}$\, 
and \,$U = \Prj^1\setminus\st{\infty} = \Afn^1$. \\ 
Let \,$\Reprf = \Autfcs{\Prj^1, 0}{\st{\infty}}$\, and  
\,$\Gp = \LalbGs{\Prj^1, 0}{\st{\infty}} 
:= \Bigpn{\Autfcs{\Prj^1, 0}{\st{\infty}}}^{\vee}$. \\ 
Using the basis \,$E_{ij} = e_i \tens e_j^{\vee} \in \gl_n(k)$\, 
we can write for $n \geq 1$: 
\begin{align*} 
	\Mor_{\Schk}\bigpn{\pn{\Afn^1,0},\pn{\gl_n, 0}} 
	&= \Biggst{\sum_{i,j,\lam} a_{ij\lam} \underbrace{E_{ij} \,t^{\lam}}_{X_{ij\lam}} \;\Bigg|\; 
		\begin{array}{l}
			1\leq i,j \leq n \\ 
			\lam \geq 1 ,\; a_{ij\lam} \in k
		\end{array}
	} 
	\laurin 
\end{align*} 
The affine algebra of \,$\LalbGs{\Prj^1, 0}{\st{\infty}}$ then becomes 
due to \Point \ref{Repr from CartierDual_Aut}, equation (\ref{O(L)_Aut}): 
\begin{align*} 
	\sO\Bigpn{\LalbGs{\Prj^1, 0}{\st{\infty}}} 
	&= \varinjlim_n \; k\lrbt{X_{ij\lam} \;\bigg|\; 
		\begin{array}{l}
			1 \leq i,j \leq n \\
			1 \leq \lam 
		\end{array}
		} \\
	&= \varinjlim_n \; k\lrbt{X_{ij\lam} \;\bigg|\; 
	\begin{array}{l}
		1 \leq i,j \leq n \\
		0 \leq \lam 
	\end{array}
	}
	\bigg/ \Bigpn{X_{ij\,0} - \del_{ij}}
\end{align*}
where $\del_{ij}$ denotes the Kronecker Delta. 
According to \cite[Remark\ 2.5.8]{Springer} \\ 
(cf.\ \Point \ref{Springer Recovery}), 
the algebra \,$\sO\bigpn{\LalbGs{\Prj^1, 0}{\st{\infty}}}$ 
is equipped with the cogroup structure
\begin{align*}
	\sO\Bigpn{\LalbGs{\Prj^1, 0}{\st{\infty}}} 
	&\lra \sO\Bigpn{\LalbGs{\Prj^1, 0}{\st{\infty}}} \tens \sO\Bigpn{\LalbGs{\Prj^1, 0}{\st{\infty}}} \\ 
	X_{ij\lam} &\lra \sum_{\substack{1\leq l\leq n \\ \mu +\nu =\lam}} 
	X_{il\nu} \tens X_{lj\mu}
\end{align*}
where because of \,$X_{ij\,0} = \del_{ij}$\, it holds 
\[ \sum_{\substack{1\leq l\leq n \\ \mu +\nu =\lam}} X_{il\nu} \tens X_{lj\mu} 
\;=\; X_{ij\lam} \tens 1 + 1 \tens X_{ij\lam} + 
\sum_{\substack{1\leq l\leq n \\ \mu +\nu =\lam;\,\mu,\nu \geq 1}} 
X_{il\nu} \tens X_{lj\mu} 
\;\laurin 
\] 
Then the geometric space with group structure is 
\begin{align*}
	\LalbGs{\Prj^1, 0}{\st{\infty}} 
	&= \Spec \sO\Bigpn{\LalbGs{\Prj^1, 0}{\st{\infty}}} \\ 
	&= \varprojlim_n \Bigpn{1 + s \gl_n[[s]]}^{\tms}
	\laurin 
\end{align*}
\eExm 

\bExm 
\label{L from HDiv for Gm}
Let \,$X = \Prj^1$, \,$S = \st{0,\infty}$\, 
and \,$U = \Prj^1\setminus\st{0,\infty} = \Afn^1 \setminus \st{0}$.  \\ 
By Example \ref{HDivf(gl) for Gm} we have a splitting of representation functors 
\[ \Autfs{\,\Prj^1, 1}{\st{0,\infty}} \cong\, \Zint \tms \Autfs{\Prj^1, 0}{\st{\infty}}  
\laurin 
\]
For \,$\LalbGs{\,\Prj^1, 1}{\st{0,\infty}} := 
\Bigpn{\Autfcs{\,\Prj^1, 1}{\st{0,\infty}}}^{\vee}$\, 
we obtain: 
\begin{align*}
	\LalbGs{\,\Prj^1, 1}{\st{0,\infty}}
	&= \Bigpn{\Autfcs{\,\Prj^1, 1}{\st{0,\infty}}}^{\vee} \\ 
	&\cong \Bigpn{\Autfcs{\Prj^1, 0}{\st{\infty}}}^{\vee} \tms \,\Zint^{\vee} \\ 
	&= \,\LalbGs{\,\Prj^1, 0}{\st{\infty}} \tms \;\Gm 
\end{align*}
using Tannaka duality (Theorem \ref{Tannaka duality}), 
Remark \ref{gbm iso}  
and Example \ref{L from HDiv for A^1}. 
\eExm

\section{Universal NonCommutative Affine Group}
\label{sec:Universal-Group}

Let $(X,\bpt)$ be a projective variety with a base point 
over a perfect field $k$. 
Let $U$ be an affine subvariety of $X$, 
and let $S := X \setminus U$ be the border. 
Although the environment is different, 
the constructions in this section are similar to the ones in \cite[Sec.\ 2]{Ru13}, 
with some simplifications due to the fact that we work with a base point 
and consider only affine groups here.

\subsection{Induced Transformation} 
\label{sub:Induced-Trafo}

Let $\Gp$ be an affine group, \emph{not} necessarily commutative, 
with unit $e$. 

\bDef 
\label{based rational maps}
A rational map \,$\phe: X \dra \Gp$\, regular at \,$\bpt$\, with \,$\phe(\bpt) = e$\, 
will be called a \emph{based rational map} from $X$ to $\Gp$, 
the set of which is denoted by 
\begin{align*} 
\Gp_{X,\bpt} 
&= \RatMap\bigpn{\pn{X,\bpt}, \pn{\Gp, e}} \\
&= \ker\Bigpn{\Gp\bigpn{\sO_{X,\bpt}} \lra \Gp\bigpn{\resfld(\bpt)}}
\end{align*} 
where \,$\resfld(\bpt)$\, is the residue field at $\bpt$. 
Then we have a natural pairing 
\begin{align*} 
\lrpair{\llul,\lull} : \Repf_{\Gp} \tms \Gp_{X,\bpt} &\lra \Autf_{X,\bpt} \\ 
\pn{\rho,\phe} &\lmt \rho \circ \phe 
\laurin 
\end{align*} 
\eDef 

\bDef 
\label{induced Trafo}
Let $\phe: X \dra \Gp$ be a based rational map to an affine group $\Gp$. 
Then \,$\trafo_{\phe}: \Gpd \lra \Autfc_{X,\bpt}$\, denotes 
the transformation induced by \,$\bigpair{\lul,\phe}$, 
where $\Gpd = \bigbt{\Repf_{\Gp}}$ is the Tannaka dual. 
For every finite dimensional $\fld$-vector space $V$ and 
\,$\lam \in \Repf_{\Gp}(V) = \Homgk\bigpn{\Gp, \Gl_V}$ we have 
\,$\trafo_{\phe}\bigpn{[\lam]} = \bt{\phe^*\lam} = \bt{\lam \circ \phe}$. 
Here $\trafo_{\phe}$ is a homomorphism of representation functors, 
due to Proposition \ref{hom of repr functors}. 
\eDef 

\bRmk 
\label{sloppy notation}
We will use representatives in order to deal with representation functors 
(cf.\ Remark \ref{G(X) for reprf}), 
and sometimes omit the brackets in the notation. 
\eRmk

\bLem 
\label{subgroup generated by X}
Let \,$\phe: X \dra \Gp$\, be a based rational map from $X$ to an affine group $\Gp$. 
Then the image of \,$\trafo_{\phe}: \Gpd \lra \Autfc_{X,\bpt}$\, 
is the representation functor corresponding to 
the subgroup of \,$\Gp$ generated by $X$. 
\eLem 

\bPf 
By construction of $\Autfc_{X,\bpt}$, 
two representations $\rho, \sig: \Gp \lra \Gl_V$ have the same image in $\Autfc_{X,\bpt}$
if and only if they coincide on the subgroup $\pair{\im \phe} \subset \Gp$ 
generated by $X$. 
Since Tannaka duality is an anti-equivalence between affine groups 
and representation functors by Theorem \ref{Tannaka duality}, 
there exist representations $\rho \neq \sig$ with the same image in $\Autfc_{X,\bpt}$ 
if and only if $\Gp \neq \pair{\im \phe}$. 
Hence the subgroup \,$\Gam \subset \Gp$\, that corresponds to the homomorphic image 
of $\Gpd$ in $\Autfc_{X,\bpt}$ satisfies \,$\Gam = \pair{\im \phe}$. 
\ePf 

\bLem
\label{Map-from-Trafo}
Let $\Gp$ be an affine algebraic group 
and \,$\trafo: \Gpd \lra \Autfc_{X,\bpt}$\, a homomorphism of representation functors. 
Then there is a unique based rational map \,$\phe: (X,\bpt) \dra (\Gp,e)$\, whose 
induced transformation is $\trafo$. 
If \,$\im(\trafo) \subset \Autfcs{X,\bpt}{\,S}$, then $\phe$ is regular outside of \,$S$, 
i.e.\ is a morphism from $U$ to $\Gp$. 
\eLem 

\bPf 
Since $\Gp$ is an affine algebraic group, it admits a faithful representation 
\,$\repr: \Gp \lra \Gl_V$ for a vector space $V$. 
Then \,$\repr \in \Homgk\bigpn{\Gp, \Gl_V} = \Gpd(V)$, 
and $\trafo(\repr) \in \Autfcs{X,\bpt}{\,S}(V) 
= \ker\bigpn{\Gl_V\bigpn{\sO_X(U)} \lra \Gl_V\pn{\resfld(\bpt)}}$ 
for some open $U \subset X$ and $S = X \setminus U$. 

Since $\rho$ is a monomorphism, there is an inverse \,$\sig: \im(\rho) \iso \Gp$\, 
such that $\sig \circ \rho = \id_{\Gp}$. 
The desired based rational map $\phe: X \dra \Gp$ 
is uniquely determined by the condition that 
$\trafo(\rho) = \trafo_{\phe}(\rho) = \rho \circ \phe$: 
indeed, $\phe  = \sig \circ \rho \circ \phe = \sig \circ \trafo(\rho)$. 
Moreover, as $\trafo$ commutes with homomorphisms of affine groups 
by Lemma \ref{hom of fmlG commutes}, we have 
\,$\trafo_{\phe}(\lam) = \lam \circ \phe = \lam \circ \sig \circ \trafo(\rho) 
= \trafo(\lam \circ \sig \circ \rho) = \trafo(\lam)$\, 
for every $\lam \in \Gpd(V)$, 
where $V$ is a finite dimensional $k$-vector space. 
\ePf 

\bLem
\label{hom of fmlG commutes}
Let $\Gp$ be an affine group 
and \,$\trafo: \Gpd \lra \Reprf$\, a homomorphism of representation functors. 
Then $\trafo$ commutes with homomorphisms of linear groups: 

If \,$\lam \in \Gpd(V) = \Homgk\pn{\Gp, \Gl_V}$\, 
factors as \,$\lam = \mu \circ \rho$\, 
for some $\rho \in \Homgk\pn{\Gp, \Gl_W}$ and 
$\mu \in \Homgk\pn{H, \Gl_V}$, 
where $V$, $W$ are finite dimensional vector spaces 
and $H = \im(\rho) \subset \Gl_W$ is a linear group, 
then \,$\trafo(\lam) = \trafo(\mu \circ \rho) = \mu \circ \trafo(\rho)$. 
\eLem 

\bPf 
By Tannaka duality, $\trafo$ admits a dual 
\,$h := \trafo^{\vee}: \Reprf^{\vee} \lra \Gp$. 
Then \,$\trafo(\mu \circ \rho) = h^{\vee}(\mu \circ \rho) = \mu \circ \rho \circ h 
= \mu \circ h^{\vee}(\rho) = \mu \circ \trafo(\rho)$. 
\ePf 


\bRmk 
\label{representative of repr subfunctor}
If $\Reprf$ is a representation subfunctor of $\Autfc_{X,\bpt}$, 
then we can find a representative \,$\Rpf \subset \Autf_{X,\bpt}$ of \,$\Reprf$. 
We will use this without further mentioning. 
\eRmk 

\bDef 
\label{support}
Let $f \in \Autf_{X,\bpt}(V)$ for some finite dimensional vector space $V$. 
Then \,$\Supp(f)$\, is the locus where $f$ is not regular. 
Let $\Rpf$ be a subfunctor of representations of $\Autf_{X,\bpt}$. Then 
\[ \Supp\pn{\Rpf} := \bigcup_{\substack{f \in \Rpf(V) \\ V \textrm{ fin.\ dim.\ vector space}}} \Supp(f) 
\laurin 
\]
\eDef

\subsubsection{Definition of a Category of Rational Maps to Torsors} 
\label{DefCatRatMaps}

\bDef 
\label{CatMr} 
A \emph{category} $\Mraff$ \emph{of rational maps from} $X$
\emph{to torsors} is a category satisfying the following conditions: 
The objects of $\Mraff$ are given by rational maps \,$\phe: X \dra P$, 
where $P$ is a torsor for an algebraic group. 
The morphisms of $\Mraff$ between two objects \,$\phe: X \dra P$\,
and \,$\psi: X \dra Q$\, are given by the set of 
those homomorphisms of torsors 
\,$h:P \ra Q$\, 
such that \,$h \circ \phe = \psi$. 
\eDef 

\bDef 
\label{Mr_H} 
Let $\Reprf$ be a representation subfunctor of $\Autfc_{X,\bpt}$. 
If $k$ is an algebraically closed field, 
then $\Mrtaff{\Reprf}{X}$ denotes the category of all 
those based rational maps $\phe: (X,\bpt) \dra (\Gp,e)$ from $X$ to affine $k$-groups 
for which the image of the induced transformation 
\,$\trafo_{\phe}:\Gpd \lra \Autfc_{X,\bpt}$
($\see$Definition \ref{induced Trafo}) lies in $\Reprf$, 
i.e.\ which induce a homomorphism of representation functors \,$\Gpd \ra \Reprf$. 

For general $k$, let $\clfld$ be an algebraic closure. 
Then we let $\Mrtaff{\Reprf}{X}$ denote the category of all 
those rational maps $\phe:X\dra P$ from $X$ to $k$-torsors 
for which the base changed map $\phe \tens \clfld$ 
is an object of $\Mrtaff{\Reprf \tens \clfld}{X \tens \clfld}$. \\ 
(Here we use 
that a torsor over $\clfld$ 
can be identified with the algebraic group acting on it, 
by means of choosing a $\clfld$-rational point.) 
\eDef 

\bDef 
\label{dual-algebraic} 
A representation functor \,$\Reprf$\, is called \emph{dual-algebraic} 
if and only if its Tannaka dual \,$\Reprfd$\, is an algebraic group. 
\eDef

\subsection{Universal Objects} 
\label{sub:Universal-Objects}

\bDef 
\label{univObj} 
Let $\Mraff$ be a category of rational maps from $X$ to torsors. 
Then \,$\lrpn{\lu:X\dra\Lu}\in\Mraff$\, is called
a \emph{universal object for} $\Mraff$ if it admits 
the universal mapping property in $\Mraff$: 
For all \,$\lrpn{\phe:X\dra P}\in\Mraff$\, there is a unique 
homomorphism of torsors \,$h:\Lu\ra P$\, such that \,$\phe=h\circ\lu$. 
\eDef 

\bRmk 
\label{univObj_unique}
Universal objects are uniquely determined up to (unique) isomorphism. 
\eRmk 

If the base field $k$ is algebraically closed, 
we can identify a torsor with the algebraic group acting on it, 
and a homomorphism of torsors becomes a homomorphism of algebraic groups 
composed with a translation (which is an isomorphism of torsors). 
In a category of based rational maps the translations are eliminated. 

\vspace{\vs} 

In the following 
we consider categories $\Mraff$ of rational maps from
$X$ to torsors 
satisfying the following condition: 

\begin{tabular}{rl}
	$\lrpn{\diamondsuit}$ & 
	If \,$\lrpn{\phe:X\dra G}\in\Mraff$\, 
	and \,$h: G \ra H$\, is a homomorphism \\ 
	& of $k$-groups, then \,$h\circ\phe \in \Mraff$\laurin 
\end{tabular} 

\bThm 
\label{Exist_univAffObj}
Assume that the base field $k$ is algebraically closed. 
Let $\Mraff$ be a category of based rational maps from $X$ to affine algebraic groups 
satisfying $\lrpn{\diamondsuit}$. 
Then the following conditions are equivalent: 

\begin{tabular}{rl}
	{\rm (i)} & For $\Mraff$ there exists a universal object 
	$\lrpn{\lu: X \dra \Lu} \in \Mraff$\laurin \\ 
	{\rm (ii)} & There is a dual-algebraic representation subfunctor $\Reprf$ of $\Autfc_{X,\bpt}$ \\ 
	& such that $\Mraff$ is equivalent to $\Mrtaff{\Reprf}{X}$\laurin \\ 
	{\rm (iii)} & The representation functor $\Reprf_{\Mraff} \subset \Autfc_{X,\bpt}$ 
	that is generated by \\ 
	& $\bigcup_{\phe \in \Mraff} \im\pn{\trafo_{\phe}}$
	is dual-algebraic and satisfies $\Mraff = \Mrtaff{\Reprf_{\Mraff}}{X}$\laurin 
\end{tabular} 

\noindent 
Here $\Mrtaff{\Reprf}{X}$ is the category of rational maps 
from $X$ to torsors under affine groups 
which induce a homomorphism of representation functors to $\Reprf$ 
($\see$Def.\ \ref{Mr_H}). 
\eThm 

\bPf  
(ii)$\Lra$(i) 
Assume that $\Mraff$ is equivalent to $\Mrtaff{\Reprf}{X}$, 
for some dual-algebraic representation subfunctor $\Reprf$ of $\Autfc_{X,\bpt}$.
In a first step we will construct an affine algebraic group $\Lu$ 
and a rational map $\lu: X \dra \Lu$. 
In a second step we will show the universality
of \,$\lu: X \dra \Lu$\, for \,$\Mrtaff{\Reprf}{X}$. 

\textbf{Step 1:} Construction of \,$\lu: X \dra \Lu \laurin$ \\
Define $\Lu$ to be the Tannaka dual of $\Reprf$, 
and the rational map \,$\lu: X \dra \Lu$\, by the condition that 
the induced transformation \,$\trafo_{\lu}: \Reprf \lra \Autfc_{X,\bpt}$\, 
from Definition \ref{induced Trafo} is the inclusion. 
This determines the map \,$\lu: X \dra \Lu$\, uniquely, 
according to Lemma \ref{Map-from-Trafo}.  
Note that \,$\lu: X \dra \Lu$\, generates $\Lu$ 
due to Lemma \ref{subgroup generated by X}. 

\textbf{Step 2:} Universality of \,$\lu: X \dra \Lu \laurin$ \\
Let \,$\phe: X \dra \Gp$\, be a rational map to an affine algebraic group 
$\Gp$, inducing a homomorphism of representation functors 
\[ \trafo_{\phe}: \Gpd \lra \Reprf \subset \Autfc_{X,\bpt}  
\hspace{15mm} \lam \lmt \lam \circ \phe 
\] 
($\see$Definition \ref{induced Trafo}). 
Let \,$h := (\trafo_{\phe})^{\vee} : \sL \ra \Gp$\, be the dual homomorphism 
of affine groups. 
We show that $\phe$ factors through $\lu$: 
\[ \xymatrix{ 
	X \ar@{-->}[r]^{\lu} \ar@{-->}@/_2pc/[rr]^{\phe} & \Lu \ar[r]^{h} & \Gp 
} 
\]
Indeed, 
\begin{eqnarray*}
	\trafo_{\,h \circ \lu} 
	& = & \pn{h \circ \lu}^* \; = \; \lu^* \circ h^{\vee} 
	\; = \; \trafo_{\lu} \circ h^{\vee} \; = \; \trafo_{\phe} 
\end{eqnarray*}
since \,$\trafo_{\lu}: \Reprf \lra \Autfc_{X,\bpt}$\, 
is the inclusion by construction of $\lu$. 
This implies that \,$h \circ \lu = \phe$\, 
by Lemma \ref{Map-from-Trafo}. 
As \,$\lu: X \ra \Lu$\, generates $\Lu$, each \,$h': \Lu \ra \Gp$\, 
fulfilling \,$h' \circ \lu = \phe$\, coincides with $h$. Hence $h$ is unique. 

(i)$\Lra$(iii) 
Assume \,$\lu: X \dra \Lu$\, is universal for $\Mraff$. 
Let $\Reprf$ be the image of the induced transformation 
\,$\trafo_\lu: \Lu^{\vee} \lra \Autfc_{X,\bpt}$. 
For $\lam\in\Lu^{\vee}(V)$  the uniqueness of the homomorphism 
\,$h_{\lam}:\Lu\ra\lam_{*}\,\Lu$\, fulfilling \,$\lam \circ \lu = h_{\lam} \circ \lu$\, 
implies that the rational maps \,$\lam \circ \lu: X \dra \lam_{*}\,\Lu$\, are 
non-isomorphic to each other for distinct \,$\lam \in \Lu^{\vee}(V)$. 
Hence \,$\trafo_{\lu}\pn{\nu} \neq \trafo_{\lu}\pn{\lam}$\, 
for \,$\nu \neq \lam \in \Lu^{\vee}(V)$. 
Therefore \,$\Lu^{\vee} \ra \Reprf$ is injective, hence an isomorphism. 

Let \,$\phe: X \dra \Gp$\, be an object of $\Mraff$. 
By universality of $\Lu$ thus  
\,$\phe: X \dra \Gp$\, factors through a unique homomorphism \,$h: \Lu \ra \Gp$. 
Then the dual homomorphism \,$h^{\vee}: \Gpd \lra \Reprf$\, 
yields a factorization of \,$\trafo_{\phe}: \Gpd \lra \Autfc_{X,\bpt}$\, 
through $\Reprf$. 
Property $\lrpn{\diamondsuit}$ 
and the existence of a universal object 
guarantee that the category $\Mraff$ contains all rational maps 
that induce a transformation to $\Reprf$, 
hence $\Mraff$ is equal to $\Mrtaff{\Reprf}{X}$.

(iii)$\Lra$(ii) 
is evident. 
\ePf 

\bNot 
\label{Lux notation}
If $\Reprf$ is a (dual-algebraic) representation subfunctor of $\Autfc_{X,\bpt}$, 
then the universal object for $\Mrtaff{\Reprf}{X}$ 
is denoted by 
\hspace{3pt}$\lalbGH{X}: X \dra \LalbGH{X}$. 

In the case $\Reprf = \Autfcs{X,\bpt}{\,S}$ 
for some closed proper subset $S \subset X$, 
we denote by $\Mrtaff{S}{X,\bpt} := \Mrtaff{\Autfcs{X,\bpt}{\,S}}{X}$ 
the category of \emph{all based morphisms} from $U = X \setminus S$ to affine groups, 
its universal object is denoted by 
\hspace{2mm}$\lalbGs{X,\bpt}{\,S}: U \lra \LalbGs{X,\bpt}{\,S}$. 

In the case $\Reprf = \Autfc_{X,\bpt}$ 
(without any specification on $\Reprf$ or $S$)
the category $\Mraff\pn{X,\bpt} := \Mrtaff{\Autfc_{X,\bpt}}{X}$ 
is the category of \emph{all based rational maps} from $X$ to affine groups, 
its universal object is denoted by 
\hspace{2mm}$\lalbG_{X,\bpt}: X \dra \LalbG_{X,\bpt}$. 
\eNot  

\bRmk 
\label{Lux construction} 
By construction, $\LalbGH{X}$ is generated by $X$. 
Since $X$ is reduced, $\LalbGH{X}$ is reduced as well, thus smooth. 
The proof of Theorem \ref{Exist_univAffObj}
showed that $\LalbGH{X}$ is the Tannaka dual of $\Reprf$. 
The rational map 
$\bigpn{\lalbGH{X}: X \dra \LalbGH{X}} \in \Mrtaff{\Reprf}{X}$ 
is characterized by the fact that 
\,$\trafo_{\,\lalbGH{X}}: \LalbGH{X}^{\vee} \lra \Autfc_{X,\bpt}$ 
is the identity $\Reprf \overset{\id} \lra \Reprf$. 
\eRmk 

\bRmk 
\label{pro maps}
Suppose $\Reprf$ is a representation subfunctor of $\Autfc_{X,\bpt}$ 
that is \emph{not} dual-algebraic, e.g.\ $\Reprf = \Autfc_{X,\bpt}$.  
Then the universal affine group $\LalbGH{X}$ is only pro-algebraic, 
and the universal map \,$\lalbGH{X}$\, 
exists only as an inverse system of rational maps 
\,$\bigpn{\lalbGa{X,\fmlF}: X \dra \LalbGa{X,\fmlF}}_{\fmlF}$, 
where $\fmlF$ runs through the dual-algebraic subfunctors of $\Reprf$ 
and $\LalbGa{X,\fmlF}$ are the algebraic quotients of $\LalbGH{X}$. 
\eRmk 

\bigskip 
Via a Galois descent we obtain over an arbitrary perfect field

\bThm 
\label{univ_affObject}
Let $\Reprf$ be a (dual-algebraic) formal $k$-subgroup of $\Autfc_{X,\bpt}$. 
The category $\Mrtaff{\Reprf}{X}$ admits a universal object 
\;$\lalbGbH{1}{X}: X \dra \LalbGbH{1}{X}$. 
Here \,$\LalbGbH{1}{X}$\, is a torsor for an (algebraic) affine group \,$\LalbGbH{0}{X}$, 
which is given by the Cartier dual of \,$\Reprf$. 
\eThm 

Functoriality of Dualization and Galois descent yields 

\bPrp 
\label{lalb_H(morphism)}
Let $\Reprf \subset \Autfc_{X,\bpt}$ be a dual-algebraic representation subfunctor. 
Let \,$\morphism: (Y,\bptt) \lra (X,\bpt)$\, be a based morphism of smooth proper varieties, 
such that $\morphism\pn{Y}$ intersects \,$\Supp(f)$ properly 
for every $f \in \Rpf(V)$, where $\Rpf \subset \Autf_{X,\bpt}$ is a representative 
and $V$ is a finite dimensional $\fld$-vector space.  
Let $\Reprff \subset \Autfc_{Y,\bptt}$ be a dual-algebraic representation subfunctor 
containing \,$\morphism^*\Reprf$. 

Then $\morphism$ induces a homomorphism of $\fld$-torsors 
\,$\LalbGba{1}{\morphism,\Reprff,\Reprf}$ 
and a homomorphism of affine algebraic $\fld$-groups 
\,$\LalbGba{0}{\morphism,\Reprff,\Reprf}$ 
\[ \LalbGba{i}{\morphism,\Reprff,\Reprf}: 
\LalbGba{i}{Y,\Reprff} \lra \LalbGba{i}{X,\Reprf} 
\hspace{8mm} \textrm{for } i = 1,0 \laurin 
\] 
\ePrp

\section{Retrieving Galois Coverings from Coverings of Affine Groups} 
\label{sec:GaloisCov-via-Isogenies}

Let $\fld$ be a finite field and $\clfld$ an algebraic closure. 

\bDef 
\label{geom Galois cover} 
A rational map $\morphism: Y \dra X$ of projective varieties over $\fld$ 
is a \emph{Galois covering} if it induces a Galois extension of function fields $K_Y|K_X$. 
It is a \emph{geometric} Galois covering if moreover 
\,$\Gal\pn{K_Y|K_X} \cong \Gal\bigpn{K_Y \clfld \,\big|\, K_X \clfld}$. 
\eDef 

\bPrp 
\label{pull-back-from-groups}
Let $(X,\bpt)$ and $(Y,\bptt)$ be pointed projective varieties over $\fld$, 
and \,$\morphism: (Y,\bptt) \dra (X,\bpt)$\, a based geometric Galois covering over $\fld$ 
with finite Galois group $\Gam$. 
Then there is a dual-algebraic representation subfunctor $\Reprf \subset \Autfc_{X,\bpt}$ 
such that the covering \,$\morphism: Y \dra X$\, is a pull-back 
of a quotient map of \,$\LalbGH{X}$\, over $\fld$, 
and \,$\Gam$ is a subgroup of \,$\LalbGH{X}\pn{\fld}$. 

Conversely, any subgroup \,$\Gam$ of a finite quotient of \,$\LalbG_{X,\bpt}\pn{\fld}$ 
defines a geometric finite Galois covering \,$\morphism: (Y,\bptt) \dra (X,\bpt)$\, 
over $\fld$. 
\ePrp 

\bPf 
If $\morphism: (Y,\bptt) \dra (X,\bpt)$ is a Galois covering with Galois group $\Gam$, 
then there are open affine subvarieties $U \ni \bpt$ of $X$ and $V \ni \bptt$ of $Y$, 
an affine algebraic group $G$ and a morphism $\phe: U \ra G/\Gam$ 
such that $V \ra U$ is isomorphic to the pull-back covering \,$\phe^*\pn{G \ra G/\Gam}$, 
see \cite[VI, No.~8, Prop.~7]{S_GrpAlg}. 
We obtain the following commutative diagram, 
see Remark \ref{pull-back construction} for an explicit construction. 
\[ \xymatrix{ 
	V \ar[r]^-{\chi} \ar[d]_-{\morphism} \ar@{}[dr]|{\square} & G \ar[d]^-{\pi}  \\ 
	U \ar[r]_-{\phe} \ar@/^0.75pc/[ur]^(0.45){\sig} & G/\Gam  
} 
\]
Then $\phe: U \dra G/\Gam$\, factors through $G$ via $\sig$, 
see Remark \ref{pull-back construction}, 
hence through a universal group $\LalbGH{X}$: 
let $\Reprf$ be the representation subfunctor of $\Autfc_{X,\bpt}$
given by the image of the induced transformation 
\,$\trafo_{\sig}: \Gd \lra \Autfc_{X,\bpt}$. 
Since $G$ is algebraic, $\Reprf$ is dual-algebraic. 
Let \,$\lu: X \dra \LalbGH{X}$\, denote the universal map. 
By the universal mapping property of $\LalbGH{X}$ there exists a homomorphism 
\,$h: \LalbGH{X} \lra G$\, such that \,$h \circ \lu = \sig$. 
Taking fibre products we obtain the following commutative diagram: 
\[ \xymatrix{ 
	V \ar[r] \ar[d]_-{\morphism} \ar@{}[dr]|{\square} & H \ar[r] \ar[d] \ar@{}[dr]|{\square} & G\tms_{G/\Gam} G \ar[r] \ar[d] \ar@{}[dr]|{\square} & G \ar[d]^-{\pi}  \\ 
	U \ar[r]_-{\lu} & \LalbGH{X} \ar[r]_-{h} & G \ar[r]_-{\pi} & G/\Gam  
} 
\] 
where \,$H = \LalbGH{X} \tms_{G/\Gam} G$. 

Now $G/\Gam$ is generated by $U$, 
see Remark \ref{pull-back construction}. 
Then Lemma \ref{generator going up} implies that $G$ is generated by $U$, 
hence a quotient of $\LalbGH{X}$. 
As algebraic quotients of $\LalbG_{X}$ correspond 
to dual-algebraic representation subfunctors of $\Autfc_{X,\bpt}$, 
there is a dual-algebraic representation subfunctor $\Reprf$ of $\Autfc_{X,\bpt}$ 
such that $G = \LalbGH{X}$. 
By construction of $G$ we have $\Gam \subset G(\fld) = \LalbGH{X}(\fld)$. 

For the converse direction, 
note that $\LalbG_{X}(\fld) = \varprojlim \LalbGH{X}(\fld)$, 
where the inverse limit ranges over 
dual-algebraic representation subfunctors $\Reprf$ of $\Autfc_{X,\bpt}$, 
and any finite quotient of $\LalbGH{X}(\fld)$ 
defines an algebraic quotient of $\LalbGH{X}$. 
Then the assertion is evident from the construction above. 
\ePf 

\bRmk[Explicit construction of pull-back of coverings] 
\label{pull-back construction} 
Let \,$\morphism: Y \dra X$\, be a Galois covering with finite group $\Gam$. 
We construct an algebraic group $G$ containing $\Gam$, 
open affines subvarieties $V \subset Y$, $U \subset X$ 
and maps $\chi$ and $\phe$ as below, 
such that $V \ra U$ is a pull-back of $G \ra G/\Gam$ 
(cf.\ \cite[VI, No.~8, Prop.~7]{S_GrpAlg}): 
\[ \xymatrix{ 
	V \ar[r]^-{\chi} \ar[d]_-{\morphism} \ar@{}[dr]|{\square} & G \ar[d]^-{\pi}  \\ 
	U \ar[r]_-{\phe} \ar@/^0.75pc/[ur]^(0.45){\sig} & G/\Gam  
} 
\]
Here \,$G = \Gm\pn{\ul{\fld}[\Gam]}$, 
with affine algebra \,$\sO(G) = \fld\bigbt{\st{X_{\gam}}_{\gam\in\Gam},\det^{-1}}$, 
where $\det := \det\bigpn{\pn{X_{\del^{-1}\gam}}_{\del,\gam\in\Gam}}$. 
There exists a normal basis $\fB$ of $K_Y$ over $K_X$, 
i.e.\ there are open subvarieties  $V \subset Y$ and $U \subset X$ 
such that $\sO(V)$ is an $\sO(U)[\Gam]$-module of rank 1, 
see \cite[VI, Thm.\ 13.1]{L_Algebra}. 
We define $\chi$ as follows: 
let $\chi^*$ map the generators of $\sO(G)$ to the normal basis $\fB$ of $K_Y\big| K_X$. 
Then $\chi$ links the algebra structure of $\sO(U)[\Gam]$ 
to the group structure of \,$G = \Gm\pn{\ul{\fld}[\Gam]}$. 
Thus $G$ is generated by $\chi: V \ra G$. 

Let $\sig = \Trace_{\morphism} \chi$ be a trace of the map $\chi$, 
defined  by $\sig(x) = \prod_{y \ra x} \chi(y)$, see \cite[III, No.~2 (b)]{S_GrpAlg}. 
(If $G$ is not commutative, the definition of $\Trace_{\morphism}\chi$ is not unique, 
but for any $x \in U$ and any fixed order of the preimages 
$y \in \morphism^{-1}(x)$ the assignment $x \lmt \prod_{y \ra x} \chi(y)$ 
can be extended locally to a well defined rational map. 
Shrinking $U$ if necessary gives a morphism \,$\sig: U \ra G$.) 
We replace $\chi$ by the $n$-fold product map $\chi^{(n)}$ defined by $y \lmt \chi(y)^n$, 
where $n = \card(\Gam)$, 
and set $\phe = \pi \circ \sig$. 
As $G/\Gam$ is generated by $G$ and the diagram commutes by construction, 
$G/\Gam$ is generated by $U$. 
\eRmk

\bLem 
\label{generator going up}
Suppose we have a pull-back diagram of Galois coverings 
as in Remark \ref{pull-back construction} above, 
with finite Galois group $\Gam$, 
where $G$ is a connected algebraic group containing $\Gam$, 
and \,$\sig = \Trace_{\morphism} \chi$\, be a trace of the map $\chi$. 

Then \,$\sig: U \ra G$\, generates \,$G$\, 
if and only if \,$\phe: U \ra G/\Gam$\, generates \,$G/\Gam$. 
\eLem 

\bPf 
Let \,$H := \pair{\im\sig}$\, be the subgroup of $G$ generated by \,$\sig: U \ra G$. 
Subgroups of algebraic groups are always closed 
(see e.g.\ \cite[1.g, Cor.\ 1.69]{Milne_algGroups}). 
If \,$\phe: U \ra G/\Gam$\, generates \,$G/\Gam$, 
then due to commutativity of the diagram $\pi(H) = G/\Gam$,  
and we have \,$\dim H \geq \dim G/\Gam = \dim G$. 
So $H$ is a closed subgroup of $G$ of the same dimension, 
and $G$ is connected, 
hence $H = G$. 
The converse direction is evident from commutativity of the diagram. 
\ePf 

\bThm 
\label{abs Galois Group}
Let $(X,\bpt)$ be a pointed projective variety over $\fld$. 
The based geometric Galois coverings of $(X,\bpt)$ 
are classified by the pro-sub-completion 
of the $\clfld$-points of the universal affine pro-algebraic group 
\,$\LalbG_{X,\bpt} = \bigpn{\Autfc_{X,\bpt}}^{\vee}$: 
\[ 
\bGalgeo{X,\bpt}
= \ProSub\bigpn{ \LalbG_{X,\bpt}\pn{\clfld} } 
\laurin 
\]
\eThm 

\bPf Due to Proposition \ref{pull-back-from-groups} 
the geometric Galois coverings of $(X,\bpt)$ that are defined over $\fld$ 
are classified by the following group: 
\begin{align}
	\bGalO{X,\bpt}
	:=& \varprojlim_{\substack{(Y,\bptt) \dra (X,\bpt) \\ \textrm{ geometric Galois}/\fld}} \Gal\pn{K_Y|K_X} \notag \\ 
	=& \,\ProSub\Biggpn{ \varprojlim_{\Reprf \subset \Autfc_{X,\bpt}} \LalbGH{X}(\fld) } 
	\label{Galois k} \mytag{G} \\ 
	=& \,\ProSub\bigpn{ \LalbG_{X}(\fld) } \laurin \notag 
\end{align} 
Here the Galois groups of geometric coverings $(Y,\bptt) \dra (X,\bpt)$ over $\fld$ 
are subgroups of $\LalbGH{X}(\fld)$ (automatically finite), 
where $\Reprf$ ranges over dual algebraic subgroups of $\Autfc_{X,\bpt}$. 
Therefore the inverse limit of these Galois groups is exactly the $\ProSub$-completion  
of $\LalbG_{X,\bpt}(\fld)$, see Notations \& Conventions \ref{pro-sub-completion}. 

The group $\bGalgeo{X,\bpt}$ 
classifies Galois coverings $(Y,\bptt)_{\clfld} \dra (X,\bpt)_{\clfld}$\, 
of the base change to the algebraic closure $\clfld$, 
i.e.\ defined over an arbitrary finite field extension $\erwfld \,|\, \fld$ in $\clfld$. 
Therefore this group is the inverse limit of the groups 
$\bGalerwO{X,\bpt}{\erwfld}$, 
where $\erwfld$ ranges over all finite field extensions $\erwfld \,|\, \fld$ in $\clfld$: 
\begin{align*}
	\bGalgeo{X,\bpt} 
	&= \varprojlim_{\erwfld | \,\fld \;\textrm{finite}} 
	   \bGalerwO{X,\bpt}{\erwfld} \\ 
	&= \varprojlim_{\erwfld | \,\fld \;\textrm{finite}} 
	   \ProSub\bigpn{ \LalbG_{X,\bpt}(\erwfld) } \\ 
	&= \,\ProSub\bigpn{ \LalbG_{X,\bpt}(\clfld) } \laurin 
\end{align*}
\hfill 
\ePf

\section{\'Etale Fundamental Group of the Affine Line}
\label{sec:FundGroup}

Instead of Galois coverings $Y \dra X$ in the sense of rational maps 
we consider now finite and \'etale Galois covers $V \ra U$ 
of an open affine subvariety $U \subset X$. 

\bQst 
In Proposition \ref{pull-back-from-groups}, 
can we specify the support of the representation functor $\Reprf$ 
(see Definition \ref{support}), 
only depending on $U$? 
Every finite dimensional $\fld$-vector space $V$ is isomorphic to $\fld^n$, 
this yields an iso \,$\Gd(V) \iso \Gd\bigpn{\fld^n} = \Homgk\pn{G,\Gl_n}$.  
For $\rho \in \Gd(V)$ it thus suffices to consider representations 
\,$\rho: G \lra \Gl_n$. 
Then \,$\rho \circ \sig: U \lra \Gl_n$\, 
is of the form \,$x \lmt \bigpn{\phe_{ij}(x)}_{ij}$\, 
for some rational functions \,$\phe_{ij} \in \sO_{X,\bpt}$\, 
such that \,$\det\bigpn{\phe_{ij}}_{ij} \in \Gm(\sO_{X,\bpt})$\, and 
\,$\trafo_{\sig}(\rho) = \lrbt{\pn{\phe_{i j}}_{i j}} 
\in \Autfc_{X,\bpt}\bigpn{\fld^n} 
= \ker\bigpn{\Gl_n(\sO_{X,\bpt}) \lra \Gl_n\pn{\resfld(\bpt)}}$. 
We have \,$\Supp(\Reprf) \subset S := X \setminus U$\, if and only if 
those rational functions satisfy 
\,$\phe_{ij} \in \sO_X(U)$\, for all $V$, $\rho \in \Gd(V)$. 
\eQst 

\bThm
\label{Abelian Coverings} 
Let $X$ be a projective curve over a finite field $\fld$. 
Let $U$ be an open affine subvariety of $X$, set $S = X \setminus U$. 
Then the abelianized geometric \'etale fundamental group of \,$U$ 
is the pro-sub-completion of the group of \,$\clfld$-points 
of the abelian universal affine pro-algebraic group 
$\Lalbs{X}{S} =\bigpn{\HDivsfc{X}{S}}^{\vee}$ (see (\ref{HDiv})): 
\[ 
\fundGabgeo{U} = \ProSub\bigpn{ \Lalbs{X}{S}\pn{\clfld} } 
\laurin 
\]
\eThm 

\bPf 
We first proceed as in Proposition \ref{pull-back-from-groups}. 
Since we consider \emph{abelian} coverings, we may use \cite[VI, No.~6, Prop.~6]{S_GrpAlg}: 
Every separable isogeny $H \ra G$ of abelian algebraic groups over a finite field 
is a quotient of the Artin-Schreier map $\asm: G \lra G$, $g \lmt \Frob(g) - g$. 
The Galois group of the covering corresponding to an isogeny is the kernel of the isogeny, 
in this case $\ker \asm = G(\fld)$, the $\fld$-valued points of $G$. 
Thus the geometric finite and \'etale abelian coverings of $U$ 
are classified by the following group: 
\begin{align*}
	\fundGabO{U} 
	:=& \varprojlim_{\substack{V \ra U \\ \textrm{ geo.\ fin.\ \'et.\ ab.\ Galois}/\fld}} \Gal\pn{K_V|K_U} \\ 
	=& \;\varprojlim_{\fmlH} \Bigpn{\ker\bigpn{\asm: \LalbH{X} \lra \LalbH{X}} } 
	= \;\varprojlim_{\fmlH} \LalbH{X}(\fld) 
\end{align*} 
where $\fmlH$ ranges over those dual-algebraic formal subgroups of $\HDivf_X$ 
that are relevant for $U$. 
Next we take care of the support: 
We need to show that the formal groups $\fmlH$ satisfy $\fmlH \subset \HDivsf{X}{S}$, 
i.e.\ the support of $\fmlH$ is contained in $S = X \setminus U$. 
For the case that $X$ is a curve 
this is proven by \cite[VI, No.~12, Prop.~11]{S_GrpAlg} and Galois descent. 
We obtain: 
\begin{align*}
	\fundGabO{U} 
	=& \;\varprojlim_{\substack{\fmlH \\ \Supp(\fmlH) \subset S}} \LalbH{X}(\fld) 
	= \,\Lalbs{X}{S}(\fld) \laurin 
\end{align*} 
Finally we consider the base change of $U$  to the algebraic closure, 
i.e.\ the field of definition is the limit over all finite extensions $\erwfld \,|\, \fld$: 
\begin{align*}
	\fundGabgeo{U} 
	=& \;\varprojlim_{\substack{\erwfld | \,\fld \\ \textrm{finite}}} \Lalbs{X}{S}(\erwfld)  
	\,=\,\ProSub\Bigpn{ \Lalbs{X}{S}(\clfld) } \laurin 
\end{align*} 
\hfill 
\ePf 

\bThm 
\label{FundGroup-of-Affine}
Let \,$X = \Prj^1$\, be the projective line over a finite field $\fld$. 
Let $\clfld$ be an algebraic closure of $\fld$. 
Let \,$U = \Afn^1$\, or \,$\Afn^1 \setminus \st{0}$, set \,$S = X \setminus U$, 
and \,$\bpt \in U$ a base point. 
Then the geometric \'etale fundamental group of $U$ 
is given by the $\clfld$-points of the universal affine pro-algebraic group \,$\LalbGs{X,\bpt}{\,S} = \bigpn{\Autfcs{X,\bpt}{\,S}}^{\vee}$: 
\[ \fundG\bigpn{U_{\clfld},\bpt} = \ProSub\bigpn{ \LalbGs{X,\bpt}{\,S}\pn{\clfld} } 
\laurin 
\]
\eThm 

\bPf 
Adapting the proof of Theorem \ref{abs Galois Group}, 
what remains to show is that in equation (\ref{Galois k}) 
for $X = \Prj^1$ and $U = \Afn^1$ or $\Afn^1 \setminus \st{0}$ 
the support of $\Reprf$ is contained in $S = X \setminus U$. 
As in the proof of Theorem \ref{Abelian Coverings} 
we can use \cite[VI, No.~12, Prop.~11]{S_GrpAlg} 
and Galois descent. 
The cited proposition only makes a statement for abelian coverings, 
but due to the explicit structure of the universal group 
\,$\LalbGs{X,\bpt}{\,S} = \bigpn{\Autfcs{X,\bpt}{\,S}}^{\vee}$\, 
from Examples \ref{L from HDiv for A^1} and \ref{L from HDiv for Gm} 
the arguments given at \textit{loc.\ cit.\ }carry over to our situation. 
We obtain: 
\begin{align}
	\fundGO{U,\bpt} 
	:=& \varprojlim_{\substack{V \ra U \\ \textrm{ geo.\ fin.\ \'et.\ Galois}/\fld}} \Gal\pn{K_V|K_U} \label{fundGO} \mytag{F} \\ 
	=& \;\ProSub\biggpn{ \varprojlim_{\substack{\Reprf \\ \Supp(\Reprf) \subset S}} \LalbGH{X}(\fld) } \notag 	
	\\ 
	=& \;\ProSub\bigpn{ \LalbGs{X,\bpt}{\,S}(\fld) } \laurin \notag
\end{align} 
Base change to the algebraic closure yields 
\begin{align*}
	\fundG\bigpn{U_{\clfld},\bpt} 
	&= \,\varprojlim_{\substack{\erwfld | \,\fld \\ \textrm{finite}}} 
	\ProSub\bigpn{ \LalbGs{X,\bpt}{\,S}(\erwfld) } \\ 
	&= \,\ProSub\bigpn{ \LalbGs{X,\bpt}{\,S}\pn{\clfld} } 
	\laurin 
\end{align*} 
\hfill 
\ePf 

\bigskip 
\bPf[Proof of Theorems \ref{FundGroup-of-A1} and \ref{FundGroup-of-Gm}] 
Follows directly from Theorem \ref{FundGroup-of-Affine}, 
taking into account the explicit description of 
\,$\LalbGs{X,\bpt}{\,S} = \bigpn{\Autfcs{X,\bpt}{\,S}}^{\vee}$\, 
for \,$X\setminus S = \Afn^1$\, from Example \ref{L from HDiv for A^1} and 
for \,$X\setminus S = \Afn^1 \setminus \st{0}$\, from Example \ref{L from HDiv for Gm}. 
\ePf

\section{Construction of Quotient Maps}
\label{sec:Quotient Maps}

Let $(X,\bpt)$ be a pointed projective variety over a finite field $\fld$, 
and $U$ an affine subvariety of $X$, set $S = X \setminus U$. 
Suppose \,$\fundG\bigpn{U_{\clfld},\bpt} 
= \ProSub\bigpn{\LalbGs{X,\bpt}{\,S}\pn{\clfld} }$. 
(This is proven for e.g.\ $U = \Afn^1$ or $\Afn^1 \setminus \st{0}$ 
by Theorem \ref{FundGroup-of-Affine}.) 

\bPnt 
\label{quotient map}
Let \,$\morphism: V \ra U$\, be a finite and \'etale Galois covering with group $\Gam$. 
We seek to describe the quotient map \,$\fundGO{U} \lsur \Gam$\, explicitly, 
where $\fundGO{U}$ is the group from (\ref{fundGO}) 
classifying geometric finite and \'etale Galois coverings of $U$ defined over $\fld$. 
Passage to the limit over finite extensions $\erwfld\,|\,\fld$ 
will give us the quotient map for the geometric \'etale fundamental group \,$\fundGgeo{U} \lsur \Gam$. 

\bigskip 
\textbf{Step 1:} \;$\LalbGs{X,\bpt}{\,S}\pn{\fld} \lsur G\pn{\fld}$ \\ 
Choose a map \,$\chi: V \lra G = \Gm\pn{\ul{k}[\Gam]}$\, 
and a local trace map \,$\phe = \Trace_{\morphism} \chi: U \ra G$\, 
as in the proof of Theorem \ref{pull-back-from-groups}, 
and a dual-algebraic representation subfunctor $\Reprf \subset \Autfc_{X,\bpt}^{\,S}$ 
such that \,$\im\pn{\trafo_{\phe}} \subset \Reprf$. 
The desired homomorphism is 
the composition of the canonical quotient map 
\,$\LalbGs{X,\bpt}{\,S} \lsur \LalbG_{\Reprf}$\,
with the unique map \,$\tha$\, obtained from 
the universal mapping property of \,$\LalbG_{\Reprf} := \LalbGs{X,\Reprf}{\,S}$: 
\[
\hspace{-10pt}
\xymatrix{
	U \ar[r]^{\phe} \ar[d]_-{\lu} & G \\ 
	\LalbG_{\Reprf} \ar[ur]_-{\tha} & 
}
\xymatrix{
	\sO(U) = \Mor\pn{U,\Afn^1} & \sO(G) = \Mor\pn{G,\Afn^1} \ar[l]_-{\phe^*} \\ 
	\sO\pn{\LalbG_{\Reprf}} \subset \varinjlim \,\Sym\pn{\Morpt\pn{U,\gl_n}} \ar[u]^-{\lu^*} 
	& \varinjlim \,\Sym\pn{\Hom\pn{G,\gl_n}} \ar[l]^-{\tha^*} \ar[u]^-{\wr} 
}
\] 
with \,$\phe^*: g \lmt g \circ \phe$, 
\hspace{5pt}$\tha^*: h \lmt h \circ \phe$, 
where the bottom line of the right diagram 
is due to \Point \ref{Repr from CartierDual_Aut} and \Point \ref{Repr from CartierDual} 
for unipotent $G$. 
Here $\tha^*$ satisfies the embedding condition $\lrpn{\bigodot}$ 
from Definition \ref{Embedding conditioon}. 
Therefore \,$\tha: \LalbG_{\Reprf} \lra G$\, is dense and hence an epimorphism 
by \cite[1.g.\ Cor.\ 1.69]{Milne_algGroups}. 

\bigskip 
Since $\Gam \subset G\pn{\fld}$, by passage to the $\ProSub$-completion 
this already yields 
\[ \fundGO{U,\bpt} = \ProSub\bigpn{\LalbGs{X,\bpt}{\,S}\pn{\fld}} 
\lsur \ProSub\bigpn{G(\fld)} 
\lsur \Gam 
\laurin 
\]
However, if $G$ is abelian, 
such a cheap trick is not needed, because 
the following direct construction will provide a homomorphism in this case: 

\bigskip 
\textbf{Step 2:} \;$G\pn{\fld} \lsur \Gam$ 
\\ 
Denote \,$H = G \tms_{G/\Gam} G$\, and let \,$\gam: H \ra G$\, 
be the pull back of the canonical projection \,$G \lra G/\Gam$, 
hence a Galois cover of group $\Gam$. 
We argue as in \cite[VI, No.~6, Prop.~6]{S_GrpAlg}: 
We have \,$\Gam \subset H\pn{\fld} = \ker\pn{\asm': H \ra H}$, 
where $\asm': x \lmt x^{-1} \Frob(x)$ is the Lang map on $H$. 
Hence $\asm'$ defines by passage to the quotient a map 
\,$\alp: H/\Gam = G \lra H$. 
Then the Lang map $\asm$ on $G$ factors as 
\,$\asm = \gam \circ \alp$, 
hence $\alp$ maps as follows: 
\;$ G\pn{\fld} = \ker\pn{\asm} = \ker\pn{\gam \circ \alp} 
\overset{\alp} \lra \ker\pn{\gam} = \Gam
\laurin 
$\, 
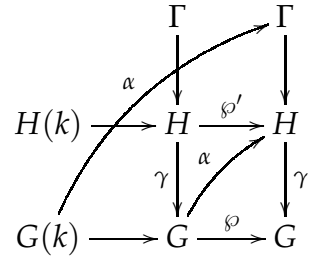
\[ \xymatrix{
	  & \Gam \ar[d] & \Gam \ar[d] \\ 
	H(\fld) \ar[r] & H \ar[r]^{\asm'} \ar[d]_-{\gam} & H \ar[d]^-{\gam} \\ 
	G(\fld) \ar@/^1.25pc/[uurr]^-{\alp} \ar[r] & G \ar@/^0.5pc/[ur]^{\alp} 
	\ar[r]^{\asm} & G
}
\] 

\bigskip 
\noindent 
If $G$ is abelian, then the Lang map $\asm$ is the Artin-Schreier map, 
which is a homomorphism. 
In this case \,$\alp: G\pn{\fld} \lsur \Gam$\, is a surjective homomorphism, 
as can be seen from the Snake-Lemma: 
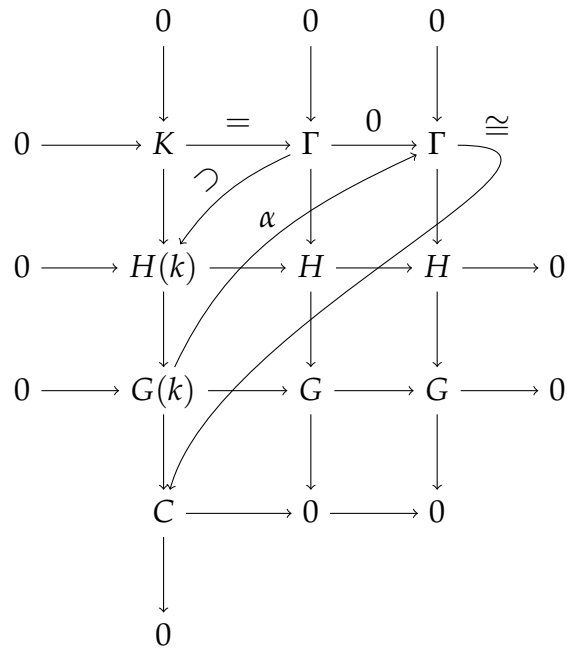
\[
\begin{tikzpicture}[text height=1.5ex, text depth=0.25ex] 
	\node (a0) at (0,0) {$0$}; 
	\node (a1) [right=of a0] {$H(\fld)$}; 
	\node (a2) [right=of a1] {$H$}; 
	\node (a3) [right=of a2] {$H$}; 
	\node (a4) [right=of a3] {$0$}; 
	\node (b0) [below=of a0] {$0$}; 
	\node (b1) [below=of a1] {$G\pn{\fld}$}; 
	\node (b2) [below=of a2] {$G$}; 
	\node (b3) [below=of a3] {$G$}; 
	\node (b4) [below=of a4] {$0$}; 
	\node (c1) [below=of b1] {$C$}; 
	\node (c2) [below=of b2] {$0$}; 
	\node (c3) [below=of b3] {$0$}; 
	\node (d1) [below=of c1] {$0$}; 
	\node (k0) [above=of a0] {$0$}; 
	\node (k1) [above=of a1] {$K$}; 
	\node (k2) [above=of a2] {$\Gam$}; 
	\node (k3) [above=of a3] {$\Gam$}; 
	\node (i1) [above=of k1] {$0$}; 
	\node (i2) [above=of k2] {$0$}; 
	\node (i3) [above=of k3] {$0$}; 
	\draw[->] 
	(i1) edge (k1)
	(i2) edge (k2)
	(i3) edge (k3)
	(k0) edge (k1) 
	(k1) edge node[above] {$=$} (k2) 
	(k2) edge node[above] {$0$} (k3) 
	(k2) edge[out=-155, in=55] node[sloped,above,pos=0.55] {$\supset$} (a1) 
	(k3) edge[out=0, in=75] node[above,pos=0.12] {$\cong$} (c1) 
	(b1) edge[out=65, in=-155] node[above,pos=0.48] {$\alp$} (k3)
	(c1) edge (c2) 
	(c2) edge (c3) 
	(a0) edge (a1) 
	(a1) edge (a2) 
	(a2) edge (a3) 
	(a3) edge (a4) 
	(b0) edge (b1) 
	(b1) edge (b2) 
	(b2) edge (b3) 
	(b3) edge (b4) 
	(k1) edge (a1) 
	(k2) edge (a2) 
	(k3) edge (a3) 
	(a1) edge (b1) 
	(a2) edge (b2) 
	(a3) edge (b3) 
	(b1) edge (c1) 
	(b2) edge (c2) 
	(b3) edge (c3) 
	(c1) edge (d1); 
\end{tikzpicture} 
\]
Since $\Gam \subset H(\fld)$, 
the map $K \ra \Gam$ is the identity, 
which implies that $\Gam \ra \Gam$ is $0$. 
Hence the snake map $\Gam \ra C$ is an isomorphism, 
because the kernel-cokernel sequence is exact. 
This implies that $\alp: G(\fld) \ra \Gam$ coincides (up to isomorphism) 
with the cokernel map $G(\fld) \ra C$, which is surjective. 

\bigskip 
Composition of the maps from Step 1 and Step 2 yields the desired projection 
\[ 
\LalbGs{X,\bpt}{\,S}(\fld) \lsur G(\fld) \lsur \Gam 
\laurin 
\]
\ePnt 

\bExm 
\label{S_3 as quotient} 
We conclude with an example for the quotient map to the Galois group of 
an unramified Galois cover of the affine line $\Afn^1$ over a finite field $\finfld_p$. 
By Theorem \ref{quasi-p-group} the Galois groups in this case are precisely 
the quasi-$p$-groups, i.e.\ groups generated by $p$-groups. 
The symmetric group $S_n$ is generated by transpositions, 
and a subgroup generated by a transposition has cardinality $2$. 
Therefore, in particular the symmetric group $S_3$ 
appears as Galois group for $\Afn^1$ over $\finfld_2$: 

Using a representation of $S_3$ by permutation matrices in $\Gl_3\pn{\finfld_2}$, 
we have an isomorphism 
\[ S_3 \iso \bigpair{T_{12}, T_{23}, T_{13}}
\]
with 
\[ 
T_{12} = \lrpn{
\begin{array}{ccc}
	0 & 1 & \\ 
	1 & 0 & \\ 
	& & 1 
\end{array}
}, \hspace{5pt}
T_{23} = \lrpn{
	\begin{array}{ccc}
		1 &  & \\ 
		& 0 & 1 \\ 
		& 1 & 0
	\end{array}
}, \hspace{5pt} 
T_{13} = \lrpn{
	\begin{array}{ccc}
		0 & & 1 \\ 
		& 1 & \\ 
		1 & & 0 
	\end{array}
}
\laurin
\]
Define \,$\Del_{\mu,\nu}$\, by the condition that 
\,$T_{\mu,\nu} = 1 + \Del_{\mu,\nu}$\,, then \,$\Del_{\mu,\nu}^2 = 0$. 
\[ 
\Del_{12} = \lrpn{
	\begin{array}{ccc}
		1 & 1 & \\ 
		1 & 1 & \\ 
		& & 0 
	\end{array}
}, \hspace{5pt}
\Del_{23} = \lrpn{
	\begin{array}{ccc}
		0 &  & \\ 
		& 1 & 1 \\ 
		& 1 & 1
	\end{array}
}, \hspace{5pt} 
\Del_{13} = \lrpn{
	\begin{array}{ccc}
		1 & & 1 \\ 
		& 0 & \\ 
		1 & & 1 
	\end{array}
}
\laurin 
\]
An \,$\finfld_2$-basis of \,$\gl_3\pn{\finfld_2}$\, is given by 
\[ \fB = \bigst{\Del_{12},\; \Del_{23},\; \Del_{13},\; 
         \Del_{12}\Del_{23},\; \Del_{23}\Del_{12},\; 
         \Del_{12}\Del_{13},\; \Del_{13}\Del_{12},\; 
         \Del_{23}\Del_{13},\; \Del_{13}\Del_{23}}
\laurin 
\]
Define the groups \,$\Gp_{\mu,\nu}$\, for \,$\st{\mu \neq \nu} \subset \st{1, 2, 3}$\, by 
\[ \Gp_{\mu,\nu}(R) = \bigst{1 + a \,\Del_{\mu,\nu} \,\big|\, a \in R} 
\] 
for every $\finfld_2$-ring $R$, 
together with morphisms 
\begin{align*}
	\phe_{\mu,\nu}: \Afn^1 &\lra \Gp_{\mu,\nu} \\ 
	a &\lmt 1 + a \,\Del_{\mu,\nu} 
\end{align*}
and let $\Gp$ be the group generated by the $\Gp_{\mu,\nu}$: 
\[ \Gp = \bigpair{\Gp_{12}, \Gp_{23}, \Gp_{13}}
\laurin 
\]
Then \,$\Gp$\, contains \,$S_3$: we have \,$\Gp_{\mu,\nu}\pn{\finfld_2} = \pair{T_{\mu,\nu}}$\, 
and \,$\Gp\pn{\finfld_2} = S_3$. 

Now consider the following algebraic quotient of \,$\LalbGs{\Prj^1, 0}{\st{\infty}}$\, 
for some \,$d \geq 1$: 
\begin{align*}
	\LalbG_{n,d} &= \Bigpn{1 + s \gl_n[[s]] \big/ \bigpn{s^{d + 1}}}^{\tms} \\ 
	\LalbG_{3,d}(R) &= \Bigst{1 + \sum_{i=1}^{d} \sum_{B\in\fB} a_i^B B \,s^i 
		\mod \bigpn{s^{d+1}} \;\Big|\; a_i^B \in R} 
	\laurin 
\end{align*}
Since the groups $\Gp_{\mu,\nu}$ are unipotent, 
by \Point \ref{quotient map} Step 1 we have the following quotient maps: 
\begin{align*}
	\tha_{\mu,\nu}: \LalbG_{3,d} &\lsur \Gp_{\mu,\nu} \\ 
	1 + \sum_{i=1}^{d} \sum_{B\in\fB} a_i^B B \,s^i 
	&\lmt 1 + a_1^{\Del_{\mu,\nu}} \Del_{\mu,\nu}
	\laurin 
\end{align*}
The composition of these maps yields a quotient map \,$\tha$\, to \,$\Gp$: 
\[ \xymatrix{ 
	\LalbG_{3,d} \ar[rr]^{\tha} \ar[rdd] \ar[rrdd] \ar[rrrdd] & & \Gp & \\ 
	 & & & \\ 
	 & \Gp_{12} \ar[uur] & \Gp_{23} \ar[uu] & \Gp_{13} \ar[uul]  
}
\]
with 
\[ \tha(g) = \tha_{12}(g) \cdot \tha_{23}(g) \cdot \tha_{13}(g) 
\]
for every \,$g \in \LalbG_{3,d}(R)$\, and every $\finfld_2$-ring $R$. 
The map $\tha$ is a priori not unique, because it depends on the order of composition. 
This is equivalent to a choice of a map 
\,$\phe = \phe_{12} \cdot \phe_{23} \cdot \phe_{13}$. 
Therefore, if we start from a morphism \,$\phe: \Afn^1 \lra G$, 
the quotient map $\tha$ is uniquely determined. 

Now for $R = \finfld_2$ we obtain the desired quotient map 
\[ \LalbG_{3,d}(\finfld_2) \lra \Gp(\finfld_2) = S_3 
\laurin 
\]
\eExm 

\bRmk 
\label{decomposition of Galois groups} 
In general, every quasi-$p$-group is a composition of $p$-groups, 
and the $p$-subgroups of a finite and \'etale group are unipotent 
(see \cite[Prop.\ 14.14]{Milne_algGroups}). 
Thus for every Galois covering 
of the affine line $\Afn^1$ over a finite field $\finfld_{p^r}$ 
the quotient map from $\LalbGs{\Prj^1, 0}{\st{\infty}}$ 
to the Galois group arises in a similar fashion as in Example \ref{S_3 as quotient}. 
\eRmk


\begin{flushright} 
	e-mail: \texttt{henrik.russell@gmail.com} 
\end{flushright}

\end{document}